\documentclass[reqno]{amsart}
\usepackage[utf8]{inputenc}
\usepackage{pgffor}
\usepackage{xcolor, soul}
\definecolor{hyperlink}{rgb}{0.7 0 0}
\usepackage[colorlinks=true,allcolors=hyperlink]{hyperref}
\usepackage{booktabs}
\usepackage{mathtools}
\usepackage{mathdots}
\usepackage{amsthm}
\usepackage{amssymb}
\usepackage{graphicx}
\usepackage{cases}
\usepackage{array}
\usepackage{enumitem}
\usepackage[figureposition=top, skip=3pt]{caption}
\usepackage{subcaption}
\usepackage[style=authoryear,
            sorting=nyt,
            backend=biber,
            natbib,
            minbibnames=7,
            maxbibnames=10]{biblatex}

\graphicspath{ {./images/} }

\usepackage[vlined,boxed,commentsnumbered]{algorithm2e} %linesnumbered

\usepackage{listings}
\usepackage[scaled]{beramono}
\usepackage[T1]{fontenc}

\usepackage{mdframed}

\newmdenv[
  topline=false,
  bottomline=false,
  rightline=false,
  skipabove=\topsep,
  skipbelow=\topsep
]{leftrule}

\newcolumntype{L}{>{$}l<{$}} % math-mode version of "l" column type
\newcommand\Tstrut{\rule{0pt}{2.6ex}}         % = `top' strut
\newcommand\Bstrut{\rule[-0.9ex]{0pt}{0pt}}   % = `bottom' strut

\renewcommand{\r}[1]{{\color{red}#1}}

\makeatletter

\DeclareCaptionSubType*{figure}
\theoremstyle{definition}
\newtheorem*{definition}{Definition}
\newtheorem{lemma}{Lemma}

\newtheorem{theorem}{Theorem}
\newtheorem{corollary}{Corollary}[theorem]

\theoremstyle{remark}
\newtheorem{remark}{Remark}[theorem]
\newtheorem*{remark*}{Remark}

\renewenvironment{proof}[1][\proofname]
    {\begin{leftrule}   
        \par        \normalfont \topsep6\p@\@plus6\p@        \trivlist
                        \item[\hskip\labelsep\bf\itshape
              #1.]\ignorespaces
    }
    {\qed\endtrivlist\end{leftrule}}
    
\renewcommand{\vec}{\mathbf}
\renewcommand*\d{\mathop{}\!\mathrm{d}}

\newcommand{\customlabel}[2]{#2\def\@currentlabel{#2}\label{#1}}

\DeclarePairedDelimiter{\norm}{|}{|}
\DeclarePairedDelimiter{\inner}{\langle}{\rangle}

\DeclareCiteCommand{\citeyear}
    {}
    {\bibhyperref{\printdate}}
    {\multicitedelim}
    {}

\DeclareCiteCommand{\citeyearpar}
    {}
    {\mkbibparens{\bibhyperref{\printdate}}}
    {\multicitedelim}
    {}
\makeatother

\calclayout

\begin{document}
\title[Tridiagonal and single-pair matrices and the inverse sum of two single-pair matrices]
    {Tridiagonal and single-pair matrices and the inverse sum of two single-pair matrices\textsuperscript{\textdagger}}

\author{S\'ebastien Bossu\textsuperscript{*}}
\thanks{*\ UNC Charlotte Department of Mathematics \& Statistics. The author thanks Stephan Sturm for helpful discussions in early versions of this paper, and an anonymous reviewer for insightful comments that have led to significant improvements.  All remaining errors are the author's responsibility.}

\date{\today\ \ \textsuperscript{\textdagger}~C++ code available at \href{https://github.com/sbossu/SPSumInverse}{github.com/sbossu/SPSumInverse}.}

\begin{abstract}
    A novel factorization for the sum of two single-pair matrices is established as product of lower-triangular, tridiagonal, and upper-triangular matrices, leading to semi-closed-form formulas for tridiagonal matrix inversion.  Subsequent factorizations are established, leading to semi-closed-form formulas for the inverse sum of two single-pair matrices. An application to derive the symbolic inverse of a particular Gram matrix is presented, and the numerical stability of the formulas is studied.
\end{abstract}



\maketitle

\section{Introduction and motivation}

For ease of notation throughout this paper, any index $i,j,k$ stated without bounds is understood to be an arbitrary natural number between 1 and $n \geq 2$.

\subsection{Single-pair matrices}      \label{sec:single-pair}

\begin{definition}
    A {\bf single-pair matrix} \citep[][pp. 78--79, 103--104]{gantmacher:2002}, also known as \emph{one-pair matrix, Green's matrix}, ``\emph{matrice factorisable}'' \citep{baranger-ducjacquet}, and \emph{generator-representable semiseparable matrix} \citep[][p.3]{vandebril:2008}, has the symmetric form
    \[
    \vec A \coloneqq \left(a_{\min(i,j)} b_{\max(i,j)} \right)
    = \begin{pmatrix}
        a_1 b_1 & a_1 b_2 & \cdots & a_1 b_n \\
        a_1 b_2 & a_2 b_2 & \cdots & a_2 b_n \\
        \vdots  & \vdots & \ddots & \vdots \\
        a_1 b_n & a_2 b_n & \cdots & a_n b_n
    \end{pmatrix}
    \]
    where $(a_i), (b_i)$ is a pair of $n$ nonzero real or complex numbers.
\end{definition}

Single-pair matrices may arise in a variety of applications, including tridiagonal matrix inversion (see \S\ref{sec:tridiag-inv} below), oscillatory mechanics, and Gram matrices --- for example:
\begin{itemize}[leftmargin=*]
    \item Covariance matrix of a Brownian sample: $\operatorname{Cov}(W_{t_i}, W_{t_j}) = t_{\min(i,j)}$ for any sampling times $0 < t_1 < \cdots < t_n$.
    \item Gram matrix of a system of step functions $u_i(x) \coloneqq 1$ if $x > x_i$, 0 otherwise, with $0 < x_1 < \cdots < x_n < 1$:
    \[
        \inner{u_i, u_j} \coloneqq \int_0^1 u_i(x) u_j(x)\d x = 1 - x_{\max(i,j)}.
    \]
\end{itemize}
Additional applications and examples may be found in \citet[ch. 3]{vandebril:2008}.

\subsubsection{The inverse of a single-pair matrix is a tridiagonal matrix} \label{sub:single-pair-inv}

The inverse of a single-pair matrix, if it exists, is known in closed form as the symmetric tridiagonal matrix \citep[p. 63]{baranger-ducjacquet}
\[
    \begin{pmatrix}
        \frac{a_2/a_1}{a_2 b_1-a_1 b_2} & \frac{-1}{a_2 b_1-a_1 b_2} 
    \\
        \frac{-1}{a_2 b_1-a_1 b_2} & \frac{a_3 b_1-a_1 b_3}{\left(a_2 b_1-a_1 b_2\right) \left(a_3 b_2-a_2 b_3\right)} & \ddots & & (0)
    \\
        & \ddots & \ddots & \ddots
    \\
        & & \ddots & \frac{a_{i+1} b_{i-1}-a_{i-1} b_{i+1}}{\left(a_i b_{i-1}-a_{i-1} b_i\right) \left(a_{i+1} b_{i}-a_{i} b_{i+1}\right)} & \ddots
    \\
        & (0) & & \ddots & \ddots & \frac{-1}{a_n b_{n-1}-a_{n-1} b_n} 
    \\
        & & & & \frac{-1}{a_n b_{n-1}-a_{n-1} b_n} & \frac{b_{n-1}/b_n}{a_n b_{n-1}-a_{n-1} b_n}
    \end{pmatrix}
\]
\begin{remark*}
    The diagonal coefficients of $\vec A^{-1}$ may be written as simple weighted sums of the surrounding off-diagonal coefficients:
    \begin{align}
        \frac{a_{i+1} b_{i-1} - a_{i-1} b_{i+1}}{\left(a_i b_{i-1}-a_{i-1} b_i\right) \left(a_{i+1} b_i - a_i b_{i+1}\right)}
        & =
        \frac{a_{i-1}/a_i}{a_i b_{i-1}-a_{i-1} b_i} + \frac{a_{i+1}/a_i}{a_{i+1} b_i - a_i b_{i+1}}
        \label{eq:sp-inv-diag-sum1}
        \\
        & = \frac{b_{i-1}/b_i}{a_i b_{i-1}-a_{i-1} b_i} + \frac{b_{i+1}/b_i}{a_{i+1} b_i - a_i b_{i+1}},
        \qquad 2 \leq i \leq n - 1.
        \label{eq:sp-inv-diag-sum2}
    \end{align}
Equation \eqref{eq:sp-inv-diag-sum2} may be extended to the cases $i = 1$ and $ i = n$ with the conventions $ a_0 = b_{n+1} = 0,\, b_0 \neq 0, a_{n+1} \neq 0$.
\end{remark*}

\subsubsection{The inverse of a symmetric tridiagonal matrix may be written as a single-pair matrix}  \label{sec:tridiag-inv}

\citet{baranger-ducjacquet} and \citet{meurant:1992} derived semi-closed-form formulas\footnote{Formulas based on a finite number of arithmetic operations, including all types of recursions.} to write the inverse of a symmetric tridiagonal matrix $\vec T \coloneqq \begin{psmallmatrix}
        \alpha_1 & -\beta_1 \\
        -\beta_1  & \alpha_2 & -\beta_2 \\
        & \ddots & \ddots & \ddots & \\
        \end{psmallmatrix} $
as a single-pair matrix
$
    \vec T^{-1} = \vec A \coloneqq \left(a_{\min(i,j)} b_{\max(i,j)} \right)
$
where the generating pair $(a_i), (b_i)$ depends on recursive coefficients parameterized by the tridiagonal coefficients of $\vec T$:
\begin{equation}    \label{eq:tridiag-inv-SP}
    \begin{dcases}
        a_i \coloneqq \frac{\beta_{i} \cdots \beta_{n-1}}{\delta_i \cdots \delta_n\,b_n}
        \\
        b_i \coloneqq \frac{\beta_1 \cdots \beta_{i-1}}{d_1 \cdots d_i}
    \end{dcases}
    ,
    \qquad
    \begin{dcases}
        d_n \coloneqq \alpha_n,\quad d_{i-1} \coloneqq \alpha_{i-1} - \frac{\beta_{i-1}^2}{d_{i}}, & i = n, n-1, \cdots, 2,
        \\
        \delta_1 \coloneqq \alpha_1, \quad
        \delta_{i+1} \coloneqq \alpha_{i+1} - \frac{\beta_{i}^2}{\delta_{i}}, & i = 1, 2, \cdots, n-1,
    \end{dcases}
\end{equation}
where any void product is deemed to be 1 \citep[pp. 710--711]{meurant:1992}.  Observe that the backward sequence $(d_i)$ must be recalculated entirely as matrix size $n$ increases.  \citet{usmani:1994} also derived semi-closed-form formulas for $\vec T^{-1}$ independently from single-pair matrices, and \citet{kilic} derived explicit inversion formulas based on backward finite continuous fractions that depend on the tridiagonal coefficients.

\subsubsection{A single-pair matrix is ``vector-bihomogeneous''}
If we multiply the pair of generators $(a_i), (b_i)$ by the same $n$ coefficients $(\lambda_i)$ elementwise, we have
\begin{equation}    \label{eq:SP-bihomogeneous}
    \operatorname{SP}( \vec {\Lambda a}, \vec{\Lambda  b} )
    = \vec\Lambda \operatorname{SP}(\vec a, \vec b)\, \vec\Lambda,
    \qquad
    \vec\Lambda \coloneqq \operatorname{diag}(\lambda_1, \cdots, \lambda_n) ,    
\end{equation}
where $\operatorname{SP}(\vec a, \vec b) \coloneqq \left(a_{\min(i,j)} b_{\max(i,j)} \right)$ denotes the single pair matrix generated by a pair of vectors $\vec a, \vec b$.

\subsection{Problem}

We are interested in calculating the inverse of the sum of two single-pair matrices in semi-closed form:
\[
    \vec {A+C} \coloneqq \left(a_{\min(i,j)} b_{\max(i,j)} + c_{\min(i,j)} d_{\max(i,j)} \right),
\]
where $a_i, b_i, c_i, d_i$ are all nonzero.  For ease of exposure, and without loss of generality, we shall set all coefficients $d_i = 1$ in the remainder of this paper, noting that the general case may be recovered by means of the factorization
\[
    \vec {A+C} = \operatorname{diag}(d_1, \cdots, d_n) \left(\frac{a_{\min(i,j)} b_{\max(i,j)}}{d_i d_j} + \frac{c_{\min(i,j)}}{d_{\min(i,j)}} \right) \operatorname{diag}(d_1, \cdots, d_n).
\]

 Closed- and semi-closed-form formulas for matrix inverses are useful for symbolic calculations, and can also help improve the speed and accuracy of numerical inversion schemes.  Sums of single-pair matrices may arise in several scenarios:
\begin{itemize}[leftmargin=*,parsep=3pt]
    \item Non-uniform perturbation $\vec C$ of an initial single-pair matrix $\vec A$;
    \item Distance matrix of an ordered set $\{ x_1 , \cdots, x_n \} $ with coefficients
    \[
        \norm{x_i - x_j} = x_{\max(i,j)} - x_{\min(i,j)};
    \]
    \item Covariance matrix for the sum of two orthogonal Brownian motions, with coefficients
        \[
            \operatorname{Cov}(W_{t_i} + Z_{\tau_i}, W_{t_j} + Z_{\tau_j}) = t_{\min(i,j)} + \tau_{\max(i,j)},
        \]
    for $t_1 < \cdots < t_n$ and $\tau_1 > \cdots > \tau_n$;
    \item Gram matrix of a system of (reverse) ReLU activation functions $ f_i(x) \coloneqq \max(0, x_i - x) $ with biases $0 < x_1 < \cdots < x_n \leq 1 $ in a single-layer, one-dimensional neural network, with coefficients (see Section \ref{sec:Gram-inv}):
    \[
        \inner{f_i, f_j} \coloneqq \int_0^1 f_i(x) f_j(x) \d x = \frac12 x_{\min(i,j)}^2 x_{\max(i,j)} - \frac16 x_{\min(i,j)}^3.
    \]
\end{itemize}

A basic solution that leverages on existing results may be obtained from the conspicuous identity
\[
    (\vec{A+C})^{-1} = \vec C^{-1} \left(\vec C^{-1} + \vec A^{-1}\right)^{-1} \vec A^{-1}
\]
where $\vec A^{-1}, \vec C^{-1}$ are symmetric tridiagonal matrices available in closed form in accordance with \S\ref{sub:single-pair-inv}.  Hence, $\vec C^{-1} + \vec A^{-1}$ is also symmetric tridiagonal and its inverse may be calculated in semi-closed form.  However, this approach has at at least two drawbacks:
\begin{itemize}[leftmargin=*,parsep=3pt]
    \item In symbolic form, the tridiagonal coefficients of $\vec C^{-1} + \vec A^{-1}$ constitute somewhat cumbersome expressions to be fed into existing tridiagonal matrix inversion formulas such as equation \eqref{eq:tridiag-inv-SP}.
    \item As observed earlier, such formulas are not incremental with matrix size $n$: all the coefficients of $\vec T^{-1}$ must be recalculated as size increases.
\end{itemize}

\subsection{Main result and organization of this paper}

Sums of single-pair matrices belong to the broader class of higher-order semiseparable matrices for which inversion algorithms have been proposed with focus on numerical applications (\cite[ch.8, 14]{vandebril:2008}, and \cite{vandebril:2007}).  Many of these existing algorithms rely on QR factorizations.  This paper derives semi-closed-form inversion formulas for the sum of two single-pair matrices satisfying certain mild conditions, based on a novel factorization.  This approach is particularly suitable for symbolic inversion, is competitive with QR for numeric inversion, and has the added benefit of being incremental with matrix size $n$.  New formulas for tridiagonal matrix inversion are also derived in the process.

The paper is organized as follows: Section 2 establishes the key factorization of this paper for the sum of two single-pair matrices as a product of lower-triangular, tridiagonal, and upper-triangular matrices, and derives corresponding formulas for tridiagonal matrix inversion.  Section 3 establishes subsequent factorizations as products of lower-triangular, single-pair, and upper-triangular matrices, and derives the main formulas for the inverse sum of two single-pair matrices (Theorems \ref{th:LSU-alt} \& \ref{th5:SP-sum-inv}).  Finally, Section 4 applies the formulas to the case of a particular Gram matrix which is how the author came to the topic of this paper, while Section 5 investigates the numerical stability of the inversion formulas.

\section{Tridiagonal decomposition of the sum of two single-pair matrices}

\begin{lemma}\label{lem:triang-prod}
    If $\mathbb L \coloneqq (l_{i,j}),\mathbb T \coloneqq (t_{i,j}), \mathbb U \coloneqq (u_{i,j})$ are respectively lower-triangular, tridiagonal and upper-triangular square matrices of order $n$, their product is an $n\times n$ matrix with coefficients
    \begin{equation}        \label{eq:LTU-product}
        (\mathbb{LTU})_{i,j} = \sum_{k=1}^{\min(i,j)} l_{i,k} t_{k,k} u_{k,j} + \sum_{k=1}^{\min(i-1,j)} l_{i,k+1} t_{k+1,k} u_{k,j} + \sum_{k=1}^{\min(i,j-1)} l_{i,k} t_{k,k+1} u_{k+1,j},
    \end{equation}
    where any void sum is deemed to be zero.  Similarly,
    \begin{equation}        \label{eq:UTL-product}
        (\mathbb{UTL})_{i,j} = \sum_{k=\max(i,j)}^n u_{i,k} t_{k,k} l_{k,j} + \sum_{k=\max(i+1,j)}^n u_{i,k-1} t_{k-1,k} l_{k,j} + \sum_{k=\max(i,j+1)}^n u_{i,k} t_{k,k-1} l_{k-1,j}.
    \end{equation}
\end{lemma}
\begin{proof} 
    By direct computation,
    \[
        (\mathbb{LT})_{i,k} = l_{i,k-1} t_{k-1, k} [2 \leq k \leq i+1] + l_{i,k} t_{k, k} [k \leq i] + l_{i,k+1} t_{k+1,k} [ k < i]
    \]
    where $[\cdot]$ denotes Iverson's bracket equal to 1 if the argument is true and 0 otherwise.  Right-multiplying by $\mathbb U$,
    \[
        (\mathbb{LTU})_{i,j} = \sum_{k=2}^{\min(i+1,j)} l_{i,k-1} t_{k-1, k} u_{k,j} + \sum_{k=1}^{\min(i,j)} l_{i,k} t_{k, k} u_{k,j} + \sum_{k=1}^{\min(i-1,j)} l_{i,k+1} t_{k+1,k} u_{k,j}.
    \]
    Rearranging and reindexing sums as appropriate yields equation \eqref{eq:LTU-product}, while similar steps yield equation \eqref{eq:UTL-product}.
\end{proof}

\begin{theorem} \label{th:single-pair-sum}
    Given a single-pair matrix $\vec A \coloneqq \left(a_{\min(i,j)} b_{\max(i,j)} \right)$  with successively distinct coefficients $b_i\neq b_{i-1}$, and a single-pair matrix $\vec C \coloneqq \left(c_{\min(i,j)} \right)$, their sum admits the factorization
    \begin{gather}
        \vec{A+C} = (\vec{LDL})\, \vec T\, (\vec{LDL})^T,
        \quad \vec L  \coloneqq \begin{pmatrix}
        1 & & (0) \\
        \vdots  & \ddots & \\
        1 & \cdots & 1
        \end{pmatrix},
        \quad \vec D  \coloneqq \operatorname{diag}(b_1 - x, b_2 - b_1, \cdots, b_n - b_{n-1}),
        \nonumber\\
        \vec T \coloneqq \begin{pmatrix}
        \alpha_1 & -\beta_1 & & & (0) \\
        -\beta_1  & \alpha_2 & -\beta_2 \\
        & \ddots & \ddots & \ddots & \\
        & & -\beta_{n-2} & \alpha_{n-1} & -\beta_{n-1} \\
        (0) & & & -\beta_{n-1} & \alpha_n
        \end{pmatrix},
        \nonumber\\
        \begin{dcases}
            \alpha_1 \coloneqq \frac{a_1 b_1 + c_1}{(b_1 - x)^2}
            \\
            \beta_1 \coloneqq \frac{a_1 x + c_1}{(b_1 - x)^2}
        \end{dcases},
        \qquad
        \begin{dcases}
            \alpha_i \coloneqq \beta_{i}+\beta_{i-1} + \frac{a_i - a_{i-1}}{b_i-b_{i-1}} - \frac{a_{i-1} - a_{i-2}}{b_{i-1}-b_{i-2}}
            \\
            \beta_i \coloneqq \frac{a_ib_{i-1} - a_{i-1}b_i + (c_i - c_{i-1})}{(b_i - b_{i-1})^2}
        \end{dcases},
        \quad 2\leq i\leq n,
        \label{eq:th1-coefs}
    \end{gather}
    where $a_0 = c_0 = 0$ and $b_0 \coloneqq x \neq b_1$ is a free parameter.
\end{theorem}
\begin{remark}
    $\vec{LDL}$ is the lower triangular matrix
     \[ \vec{LDL} =
        \begin{pmatrix}
            b_1 - x &  &  &  (0)  \\
            b_2 - x & b_2 - b_1 \\
            \vdots  & \vdots & \ddots \\
            b_n - x & b_n - b_1 & \cdots & b_n - b_{n-1}
        \end{pmatrix}
    \]    
\end{remark}
\begin{remark}
    Some natural choices for $x$ include: $x=0,\ x=b_n$ if $b_n \neq b_1$, and $x = 1-b_1$ if $b_1 \neq \frac12$.
\end{remark}
\begin{proof}
    Since $(\vec{LDL})\,\vec T\,(\vec{LDL})^T$ is clearly symmetric we may merely identify upper triangular coefficients. By  Lemma \ref{lem:triang-prod}, %, equation \eqref{eq:LTU-product},
    \begin{multline*}
        \left((\vec{LDL})\,\vec T\,(\vec{LDL})^T\right)_{i,j}
        = \sum_{k=1}^{\min(i,j)} (b_i-b_{k-1})\alpha_k (b_j - b_{k-1})
        - \sum_{k=1}^{\min(i-1,j)} (b_i-b_k)\beta_k (b_j - b_{k-1}) 
        \\
        - \sum_{k=2}^{\min(i+1,j)} (b_i-b_{k-2})\beta_{k-1} (b_j - b_{k-1}),
    \end{multline*}
    where any void sum is deemed to be zero.  We now proceed by induction: it is easy to verify that $\left((\vec{LDL})\,\vec T\,(\vec{LDL})^T\right)_{1,1} = a_1 b_1 + c_1$; if $ \left((\vec{LDL})\,\vec T\,(\vec{LDL})^T\right)_{i,j} = a_i b_j + c_i $ holds for all $i \leq j$, then
    \begin{multline*}
        \left((\vec{LDL})\,\vec T\,(\vec{LDL})^T\right)_{i,j+1}
        = 
        \left((\vec{LDL})\,\vec T\,(\vec{LDL})^T\right)_{i,j}
        \\
        + (b_{j+1}-b_j)\Bigg(\sum_{k=1}^{i} (b_i-b_{k-1})\alpha_k 
        - \sum_{k=1}^{i-1} (b_i-b_k)\beta_k 
        - \sum_{k=2}^{i+1} (b_i-b_{k-2})\beta_{k-1},
        \Bigg),
        \qquad i \leq j.
    \end{multline*}
    Reindexing sums and simplifying, then substituting equation \eqref{eq:th1-coefs} for $\alpha_k$, we have for all $ i \leq j $
    \begin{align*}
        &\sum_{k=1}^{i} (b_i-b_{k-1})\alpha_k - \sum_{k=1}^{i-1} (b_i-b_k)\beta_k 
        - \sum_{k=2}^{i+1} (b_i-b_{k-2})\beta_{k-1}
        \\
        =\ & (b_i - x) (\alpha_1 - \beta_1) + \sum_{k=2}^{i} (b_i-b_{k-1})(\alpha_k - \beta_{k-1} - \beta_k)
        \\
        =\ & (b_i - x) (\alpha_1 - \beta_1) + \sum_{k=2}^{i} (b_i-b_{k-1})\left( \frac{a_k - a_{k-1}}{b_k-b_{k-1}} - \frac{a_{k-1} - a_{k-2}}{b_{k-1}-b_{k-2}} \right)
        \\
        =\ & (b_i - x) (\alpha_1 - \beta_1) + \left( (b_i - b_{i-1})\frac{a_i - a_{i-1}}{b_i-b_{i-1}} - (b_i - b_1)\frac{a_1}{b_1 - x} \right) + \sum_{k=2}^{i-1} (b_k-b_{k-1})\frac{a_k - a_{k-1}}{b_k-b_{k-1}}
    \end{align*}
    where we applied summation by parts (Abel's lemma) in the last step.  Substituting the definitions of $\alpha_1, \beta_1$, the above expression simplifies to $a_i$ and we recover
    \begin{align*}
        \left((\vec{LDL})\,\vec T\,(\vec{LDL})^T\right)_{i,j+1}
        & =\  
        \left((\vec{LDL})\,\vec T\,(\vec{LDL})^T\right)_{i,j} + (b_{j+1} - b_j)a_i
        \\
        & =\ a_i b_j + c_i + (b_{j+1} - b_j)a_i = a_i b_{j+1} + c_i,
        \qquad i \leq j,
    \end{align*}
    as required.  Finally, for $i = j+1$,
    \begin{multline*}
        \left((\vec{LDL})\,\vec T\,(\vec{LDL})^T\right)_{j+1,j+1}
        = \sum_{k=1}^{j+1} (b_{j+1}-b_{k-1})\alpha_k (b_{j+1}-b_{k-1})
        - \sum_{k=1}^{j} (b_{j+1}-b_k)\beta_k (b_{j+1}-b_{k-1})
        \\
        - \sum_{k=2}^{j+1} (b_{j+1}-b_{k-2})\beta_{k-1}(b_{j+1}-b_{k-1}).
    \end{multline*}
Substituting $b_{j+1} = (b_{j+1} - b_j) + b_j$ inside each first factor, distributing, rearranging, and attending to upper and lower summation bounds,
    \begin{flalign*}
        & \mathrlap{ \left((\vec{LDL})\,\vec T\,(\vec{LDL})^T\right)_{j+1,j+1} }
        \\
        =&\ \mathrlap{\sum_{k=1}^{j} (b_{j}-b_{k-1})\alpha_k (b_{j+1}-b_{k-1})
        - \sum_{k=1}^{j-1} (b_{j}-b_k)\beta_k(b_{j+1}-b_{k-1})
        - \sum_{k=2}^{j+1} (b_{j}-b_{k-2})\beta_{k-1}(b_{j+1}-b_{k-1})}
        \\
        && + (b_{j+1}-b_j)\Bigg( (b_{j+1}-x)(\alpha_1 - \beta_1)
        + (b_{j+1}-b_j)\beta_{j+1} + 
        \sum_{k=2}^{j+1} (b_{j+1}-b_{k-1})(\alpha_k - \beta_k 
        - \beta_{k-1})\Bigg).
    \end{flalign*}
    Recognizing $ \left((\vec{LDL})\,\vec T\,(\vec{LDL})^T\right)_{j,j+1}$, following the same steps as before to simplify the sum between brackets, and simplifying,
    \begin{align*}
        \left((\vec{LDL})\,\vec T\,(\vec{LDL})^T\right)_{j+1,j+1}
        & = \left((\vec{LDL})\,\vec T\,(\vec{LDL})^T\right)_{j,j+1}
        + (b_{j+1}-b_j)a_j + (b_{j+1}-b_j)^2\beta_{j+1}
        \\
        & = a_{j+1}b_{j+1} + c_{j+1},
    \end{align*}
    where we substituted the definition of $\beta_{j+1}$ and simplified in the last step.
\end{proof}

\begin{corollary}
    $
        \det(\vec{A+C}) = (b_1 - x)^2 (b_2 - b_1)^2 \cdots (b_n - b_{n-1})^2 \det\vec T
    $, and $\vec{A+C}$ is invertible if and only if $\vec T$ is invertible.
\end{corollary}

\noindent The case when $\vec C$ is the null matrix (i.e. $ c_i = 0$) leads to a new factorization for symmetric tridiagonal matrices:
\begin{theorem}     \label{th:tridiagonal-single-pair}
    Subject to the generalized continuant\footnote{A \emph{generalized continuant} sequence $(b_i)$ obeys a three-term recurrence relation of the form $ b_i = \nu_i b_{i-1} - \omega_{i-1} \rho_{i-1} b_{i-2} $ \citep[pp. 516--525]{muir:1960}.  In this instance, we have $\nu_i = \alpha_i / \beta_i, \omega_i = \beta_i, \rho_i = 1/\beta_{i+1}$.} sequence $(b_i)$ defined below satisfying $b_i \neq b_{i-1}$, an irreducible\footnote{A symmetric tridiagonal matrix is irreducible when all its off-diagonal coefficients are nonzero.} symmetric tridiagonal matrix $\vec T \coloneqq \begin{psmallmatrix}
        \alpha_1 & -\beta_1 \\
        -\beta_1  & \alpha_2 & -\beta_2 \\
        & \ddots & \ddots & \ddots & \\
        \end{psmallmatrix} $
    admits the factorization
    \[
        \vec T = (\vec{LDL})^{-1}\, \vec A\, (\vec{LDL})^{-T},
        \quad
        \vec A \coloneqq \left(a_{\min(i,j)} b_{\max(i,j)} \right),
        \quad 
        \vec D  \coloneqq \operatorname{diag}(b_1 - x, b_2 - b_1, \cdots, b_n - b_{n-1}),
    \]
    with the generating pair $(a_i), (b_i)$ determined by the second order coupled recursion
    \begin{equation}    \label{eq:th2-recursion-init}
        \begin{dcases}
            a_1 \coloneqq x\frac{(\alpha_1 - \beta _1 )^2}{\beta_1}
            \\
            b_1 \coloneqq x\frac{\alpha_1}{\beta_1}
        \end{dcases},
    \end{equation}
    and, for $2 \leq i \leq n $,
    \begin{numcases}{}
        a_i \coloneqq \ a_{i-1}  -  \frac {a_{i-1}}{\beta_i} \sum_{j=1}^i(\alpha _j - \beta _j - \beta _{j-1})
        + \frac{b_{i-1}}{\beta_i}\left(\sum_{j=1}^i(\alpha _j - \beta _j - \beta _{j-1}) \right)^2
        \label{eq:th2-recursion-a}
        \\
        b_i \coloneqq \ \frac{\alpha_i}{\beta_i} b_{i-1}  - \frac{\beta_{i-1}}{\beta_i}b_{i-2},
        \label{eq:th2-recursion-b}
    \end{numcases}
    where $a_0 = \beta_0 = 0$, and $b_0 \coloneqq x \neq 0, \; \beta_n \coloneqq y \neq 0$ are free parameters.
\end{theorem}
\begin{remark}
    If $\alpha_i = \beta_i + \beta_{i-1}$ for all $i$ then $\vec T$ is singular with $(1, \cdots, 1)^T \in \ker \vec T$, and $a_i = a_{i-1}, b_i = b_{i-1}$ for all $i$.
\end{remark}
\begin{remark}
    Under this factorization, the calculation of the generating pair $(a_i), (b_i)$ is incremental with matrix size $n$.
\end{remark}
\begin{remark}
    $(\vec{LDL})^{-1}$ is the lower triband matrix
    \[  (\vec{LDL})^{-1} =
        \begin{pmatrix}
            \frac{1}{b_1-x} &  &  &  (0)  \\
         -\frac{b_2-x}{\left(b_1-x\right) \left(b_2-b_1\right)} & \frac{1}{b_2-b_1} 
         \\
          \frac{1}{b_2-b_1} & -\frac{b_3-b_1}{\left(b_2-b_1\right) \left(b_3-b_2\right)} & \frac{1}{b_3-b_2}
          \\
          (0) & \ddots & \ddots & \ddots
        \end{pmatrix}
        = \begin{bmatrix}
            \vec I & \vec 0
        \end{bmatrix} \left(b_{\max(i,j)-1}\right)_{1\leq i,j\leq n+1}^{-1} \begin{bmatrix}
            \vec 0^T \\ \vec I
        \end{bmatrix}
    \]    
    where $\vec I\in\mathbb R^{n\times n}$ is the identity matrix and $\vec 0\in\mathbb R^n$ is the zero vector.  As such, $(\vec{LDL})^{-1}$ may be viewed as a shifted inverse of a particular single-pair matrix.
\end{remark}

\begin{proof}[Proof of Theorem \ref{th:tridiagonal-single-pair}]
      Having established the factorization of Theorem \ref{th:single-pair-sum}, we only need to show that the coefficients $ a_i, b_i $ recursively defined by equations \eqref{eq:th2-recursion-a} and \eqref{eq:th2-recursion-b} solve equation \eqref{eq:th1-coefs} when $\vec C = \vec O$.  We proceed by induction: it is easy to see that $a_1, b_1$ satisfy equation \eqref{eq:th2-recursion-init}.  Suppose that equation \eqref{eq:th1-coefs} is satisfied up to index $i-1$; with some algebra, it can be shown that solving equation~\eqref{eq:th1-coefs} for $a_i, b_i$ as unknowns yields
    \begin{numcases}{}
            \hat a_i  \coloneqq \ a_{i-1}  -  \frac {a_{i-1}}{\beta_i} \left(\frac{a_{i-1}-a_{i-2}}{b_{i-1}-b_{i-2}} + \alpha _i - \beta _i - \beta _{i-1} \right)
            + \frac{b_{i-1}}{\beta_i}\left(\frac{a_{i-1}-a_{i-2}}{b_{i-1}-b_{i-2}} + \alpha _i - \beta _i - \beta _{i-1} \right)^2
            \label{eq:th2-proof-rec-a}
            \\
            \hat b_i \coloneqq \ b_{i-1}\left(1 + \frac{\alpha_i -\beta_i - \beta_{i-1}}{\beta_i} \right)  +  \frac1{\beta_i}\frac{a_{i-1} b_{i-2}-a_{i-2} b_{i-1}}{b_{i-1}-b_{i-2}}
            \label{eq:th2-proof-rec-b}
    \end{numcases}
    \begin{itemize}[leftmargin=1em]
        \item By induction assumption, we have $ \displaystyle \beta_{i-1} = \frac1{\beta_i}\frac{a_{i-1} b_{i-2}-a_{i-2} b_{i-1}}{\left(b_{i-1}-b_{i-2}\right)^2} $.  Substituting into \eqref{eq:th2-proof-rec-b},
        \[
            \hat b_i = \ b_{i-1}\left(1 + \frac{\alpha_i -\beta_i - \beta_{i-1}}{\beta_i} \right)  +  \frac{\beta_{i-1}}{\beta_i}(b_{i-1}-b_{i-2}),
        \]
        which simplifies to $ \frac{\alpha_i}{\beta_i} b_{i-1}  - \frac{\beta_{i-1}}{\beta_i}b_{i-2} =: b_i $ as required.
                \item By induction assumption for all $ 2 \leq j \leq i - 1$,
        \[
            \alpha_j = \beta_{j-1}+\beta_j + \frac{a_j - a_{j-1}}{b_j-b_{j-1}} - \frac{a_{j-1} - a_{j-2}}{b_{j-1}-b_{j-2}},
            \qquad 2 \leq j \leq i - 1,
        \]
        and for $ j = 1 $ we have $ \alpha_1 - \beta_1 - \beta_0 = \frac{a_1 - a_0}{b_1-b_0} $ with the convention $ a_0 = \beta_0 = 0$ and $ b_0 = x $.  Thus,
        \[
            \sum_{j=1}^{i-1} (\alpha_j - \beta_j - \beta_{j-1}) = \frac{a_{i-1} - a_{i-2}}{b_{i-1}-b_{i-2}}.
        \]
        Substituting into equation \eqref{eq:th2-proof-rec-a},
        \begin{align*}
            \hat a_i  = \ a_{i-1}  & -  \frac {a_{i-1}}{\beta_i} \left(\sum_{j=1}^{i-1} (\alpha_j - \beta_j - \beta_{j-1}) + \alpha _i - \beta _i - \beta _{i-1} \right)
            \\
            & + \frac{b_{i-1}}{\beta_i}\left(\sum_{j=1}^{i-1} (\alpha_j - \beta_j - \beta_{j-1}) + \alpha _i - \beta _i - \beta _{i-1} \right)^2,
        \end{align*}
        which matches equation \eqref{eq:th2-recursion-a} as required.
    \end{itemize}
    \vspace{-20pt}
\end{proof}

\begin{corollary}   \label{corol:tridiag-inverse}
    If $(b_i)$ satisfies $b_i \notin \{0, b_{i-1}\}$, then $\vec T$ is invertible in semi-closed form as $ \vec T^{-1} = (\vec{LDL})^T\, \vec A^{-1}\, (\vec{LDL}) $, where $\vec A^{-1}$ is the symmetric tridiagonal matrix given in subsection \ref{sec:single-pair},
     \[
        \vec A^{-1} = \begin{pmatrix}
        \lambda_1 & -\mu_1 & & & (0) \\
        -\mu_1  & \lambda_2 & -\mu_2 \\
        & \ddots & \ddots & \ddots & \\
        & & -\mu_{n-2} & \lambda_{n-1} & -\mu_{n-1} \\
        (0) & & & -\mu_{n-1} & \lambda_n
        \end{pmatrix},
    \]
    with tridiagonal coefficients
    \begin{numcases}{}
        \lambda_i \coloneqq
        \frac{b_{i+1}}{b_i}\mu_{i} + \frac{b_{i-1}}{b_i}\mu_{i-1},
        & $ 1 \leq i \leq n,$
        \label{eq:single-pair-inverse-lambda}
        \\
        \mu_i \coloneqq \frac1{a_{i+1} b_{i} - a_{i} b_{i+1}}
        = \frac1{(b_{i+1}-b_{i})^2\beta_{i+1}},
        & $ 1 \leq i \leq n - 1,$
        \label{eq:single-pair-inverse-mu}
                            \end{numcases}
    where $\mu_0 \coloneqq \frac1{a_1 x} = \frac{\beta_1}{x^2(\alpha_1 - \beta_1)^2},\, \mu_n \coloneqq 0$, and $b_0 \coloneqq x \neq 0, \; \beta_n \coloneqq y \neq 0$ are free parameters.  \end{corollary}
\begin{remark}
    The sequence $(a_i)$ is not required to calculate the coefficients of $\vec T^{-1}$.  This apparent ``loss of information symmetry'' may be disconcerting at first sight: if $\vec T = (\vec{LDL})^{-1}\, \vec A\, (\vec{LDL})^{-T}$ where the single-pair matrix $\vec A$ is generated by the pair $(a_i), (b_i)$ of $2n$ recursive coefficients that depend on $\alpha$'s and $\beta$'s, one might expect $\vec T^{-1}$ to be equally generated by the same $2n$ recursive coefficients.  It turns out that we may construe the sequence $a_1, \dots, a_n$ determined by equation \eqref{eq:th2-recursion-a} as a first-order non-homogeneous recurrence with variable coefficients that depend on $b$'s, $\alpha$'s and $\beta$'s:
    \[
        a_i =  \underbrace{\left[1 -  \frac1{\beta_i} \sum_{j=1}^i(\alpha _j - \beta _j - \beta _{j-1})\right]}_{f_i} a_{i-1}
        + \underbrace{\left[\frac1{\beta_i}\left(\sum_{j=1}^i(\alpha _j - \beta _j - \beta _{j-1}) \right)^2\right] b_{i-1}}_{g_i} 
    \]
with known solution \citep[p. 46]{mickens:2015}
    \[
        a_i = a_1 \prod_{j=1}^{i-1} f_j + \sum_{j=1}^{i-1} g_j \prod_{k=j+1}^{i-1} f_k^{-1}.
    \]
    In other words, the sequence $a_1, \dots, a_n$ is subsidiary to $b_1, \dots, b_n$.
\end{remark}
\begin{corollary}       \label{corol:tridiagonal-single-pair-2}
    At the expense of slightly more complex surrounding factors, we can take advantage of the ``vector bihomogeneity'' of single-pair matrices (equation \ref{eq:SP-bihomogeneous}) to obtain the alternative factorization
    \begin{gather*}
        \vec T = \big( \vec L^{-1} \vec{\dot D}^{-1} \vec B \big)\: \vec{\dot A}\: \big( \vec L^{-1} \vec{\dot D}^{-1} \vec B \big)^T,
        \qquad \vec{\dot A} \coloneqq \left(\dot a_{\min(i,j)} \dot b_{\max(i,j)} \right),
        \quad \vec{\dot D} \coloneqq \operatorname{diag}( \dot b_i - \beta_i \dot b_{i-1} ),
        \\
        \vec B \coloneqq  \begin{pmatrix}
        1 & & &  (0) \\
        - \beta_2 & 1  \\
        & \ddots & \ddots  \\
        (0) &  & -\beta_n & 1
        \end{pmatrix},
        \qquad \beta_n \coloneqq y \neq 0,
    \end{gather*}
    with the generating pair $(\dot a_i), (\dot b_i)$ determined by the second order coupled recursion
    \[
        \begin{dcases}
            \dot a_1 \coloneqq x(\alpha_1 - \beta _1 )^2
            \\
            \dot b_1 \coloneqq x\: \alpha_1
        \end{dcases},
    \]    
    and, for $2 \leq i \leq n $,
    \begin{numcases}{}
        \dot a_i \coloneqq \ \beta_i \dot a_{i-1}  -  \dot a_{i-1} \sum_{j=1}^i(\alpha _j - \beta _j - \beta _{j-1})
        + \dot b_{i-1}\left(\sum_{j=1}^i(\alpha _j - \beta _j - \beta _{j-1}) \right)^2
        \label{eq:corol2.2-recursion-a}
        \\
        \dot b_i \coloneqq \ \alpha_i \dot b_{i-1}  - \beta_{i-1}^2 \dot b_{i-2},
        \label{eq:corol2.2-recursion-b}
    \end{numcases}
    where $a_0 = \beta_0 = 0$, and $b_0 \coloneqq x \neq 0, \; \beta_n \coloneqq y \neq 0$ are free parameters.
\end{corollary}
\begin{remark}
    $\dot b_i$ is the generalized continuant associated with the symmetric tridiagonal coefficients $\alpha_1, \cdots, \alpha_i, -\beta_1, \cdots, -\beta_{i-1}$, i.e. it is the determinant of the leading principal submatrix of order $i$ of $\vec T$.
\end{remark}
\begin{proof}[Proof of Corollary \ref{corol:tridiagonal-single-pair-2}]
    Substitute equation $\vec A = \vec\Lambda^{-1} \vec{\dot A} \vec\Lambda^{-1},\ \vec\Lambda \coloneqq \operatorname{diag} \left(\prod_{j=1}^i \beta_j\right) $ into the factorization of Theorem \ref{th:LSU}, and simplify.
\end{proof}

\section{Single-pair decomposition of the sum of two single-pair matrices}

\begin{theorem}     \label{th:LSU}
    Subject to the generalized continuant sequence $(v_i)$ defined below satisfying $v_i \neq v_{i-1}$, the sum of two single-pair matrices $\vec A \coloneqq \left(a_{\min(i,j)} b_{\max(i,j)} \right), \vec C \coloneqq \left(c_{\min(i,j)} \right)$ with $b_i\neq b_{i-1}$, admits the factorization
    \[
        \vec{A+C} = (\vec{L\Delta L}^{-1})\, \vec M\, (\vec{L\Delta L}^{-1})^T,
        \qquad \vec \Delta \coloneqq \operatorname{diag} \left( \frac{b_i - b_{i-1}}{v_i-v_{i-1}} \right),
        \qquad \vec M \coloneqq \left(u_{\min(i,j)} v_{\max(i,j)} \right),
    \]
    with the generating pair $(u_i), (v_i)$ determined by the second-order coupled recursion
    \begin{equation}    \label{eq:th3-coefs-init}
        \left\{\begin{aligned}            u_1 & \coloneqq z\frac{a_1^2}{a_1 x + c_1}
            \\
            v_1 & \coloneqq z\frac{a_1b_1 + c_1}{a_1 x + c_1}
        \end{aligned}\right. ,        
    \end{equation}
    and, for $2 \leq i \leq n $,
    \begin{numcases}{}
            u_i \coloneqq
             u_{i-1} - \frac{(b_i - b_{i-1})^2}{a_i b_{i-1}-a_{i-1} b_i +c_i-c_{i-1}} \left[
                u_{i-1} \frac{a_i - a_{i-1}}{b_i-b_{i-1}}
                - v_{i-1} \left( \frac{a_i - a_{i-1}}{b_i-b_{i-1}} \right)^2
                \right]
            \label{eq:th3-coefs-u}
            \\
            v_i \coloneqq v_{i-1} + \frac{(b_i - b_{i-1})^2}{a_i b_{i-1}-a_{i-1} b_i +c_i-c_{i-1}}
            \nonumber\\ 
            \qquad\quad \times \left[ v_{i-1}\left( \frac{a_i - a_{i-1}}{b_i-b_{i-1}} - \frac{a_{i-1} - a_{i-2}}{b_{i-1}-b_{i-2}} \right)
            + (v_{i-1}-v_{i-2}) \frac{a_{i-1} b_{i-2}-a_{i-2} b_{i-1} +c_{i-1}-c_{i-2} }{(b_{i-1}-b_{i-2})^2} \right],
            \label{eq:th3-coefs-v}
    \end{numcases}
    where $a_0 \coloneqq c_0 \coloneqq 0$, and $b_0 \coloneqq x \neq b_1, v_0 \coloneqq z \neq 0$ are free parameters such that $a_1 x + c_1 \neq 0$.
\end{theorem}
\begin{remark}
    $\vec{L\Delta L}^{-1}$ is the lower triangular matrix
    \[
        \vec{L\Delta L}^{-1} = \begin{pmatrix}
        \delta_1 \\
        \delta_1 - \delta_2  & \delta_2 & & (0)  \\
        \delta_1 - \delta_2  & \delta_2 - \delta_3 & \delta_3  \\
        \vdots & \vdots & \ddots & \ddots \\
        \delta_1 - \delta_2 & \delta_2 - \delta_3 & \cdots & \delta_{n-1} - \delta_n  & \delta_n
        \end{pmatrix},
        \qquad \delta_i \coloneqq \frac{b_i - b_{i-1}}{v_i-v_{i-1}}.
    \]
\end{remark}
\begin{proof}
    From Theorem \ref{th:single-pair-sum}, $ \vec{A+C} = (\vec{LD_1 L})\, \vec T\, (\vec{LD_1 L})^T, 
    \ \vec D_1  \coloneqq \operatorname{diag}(b_1 - x, b_2 - b_1, \cdots, b_n - b_{n-1}) $.  From Theorem \ref{th:tridiagonal-single-pair},
    $
        \vec T = (\vec{LD_2 L})^{-1}\, \vec M\, (\vec{LD_2 L})^{-T},
        \ \vec D_2  \coloneqq \operatorname{diag}(v_1 - z, v_2 - v_1, \cdots, v_n - v_{n-1}),
        \ \vec M \coloneqq \left(u_{\min(i,j)} v_{\max(i,j)} \right)
    $.  Combining and simplifying yields the desired factorization, with
    \[
        \begin{dcases}
            u_1 \coloneqq z\frac{(\alpha_1 - \beta _1 )^2}{\beta_1}
            \\
            v_1 \coloneqq z\frac{\alpha_1}{\beta_1}
        \end{dcases},
    \]    
    and, for $2 \leq i \leq n $,
    \begin{numcases}{}
        u_i \coloneqq \ u_{i-1}  -  \frac {u_{i-1}}{\beta_i} \left(\sum_{j=1}^i (\alpha_j - \beta_j - \beta_{j-1}) \right)
            + \frac{v_{i-1}}{\beta_i}\left(\sum_{j=1}^i (\alpha_j - \beta_j - \beta_{j-1}) \right)^2
        \label{eq:th3-proof-recursion-u}\\
        v_i \coloneqq \ \frac{\alpha_i}{\beta_i} v_{i-1} -  \frac{\beta_{i-1}}{\beta_i}{v_{i-2}}
        \label{eq:th3-proof-recursion-v}
    \end{numcases}
    Substituting equation \eqref{eq:th1-coefs} for $\alpha_i, \beta_i, \beta_{i-1}$ into the above expressions, we recover equations \eqref{eq:th3-coefs-init}, \eqref{eq:th3-coefs-u} and \eqref{eq:th3-coefs-v} after some algebra.
\end{proof}

\begin{corollary}       \label{corol:sum-inv}
    If $(v_i)$ satisfies $v_i \notin \{0, v_{i-1}\}$, then $\vec{A+C}$ is invertible in semi-closed form as
    \[
        (\vec{A+C})^{-1} = \left(\vec{L\Delta}^{-1}\vec L^{-1}\right)^T \vec M^{-1} \left(\vec{L\Delta }^{-1}\vec L^{-1}\right),
        \qquad
        \vec M^{-1} =: \begin{psmallmatrix}
        \lambda_1 & -\mu_1 & & & (0) \\
        -\mu_1  & \lambda_2 & -\mu_2 \\
        & \ddots & \ddots & \ddots & \\
        & & -\mu_{n-2} & \lambda_{n-1} & -\mu_{n-1} \\
        (0) & & & -\mu_{n-1} & \lambda_n
        \end{psmallmatrix},
    \]
    with tridiagonal coefficients
    \begin{numcases}{}
        \lambda_i \coloneqq
         \frac{v_{i+1}}{v_i}\mu_{i} + \frac{v_{i-1}}{v_i}\mu_{i-1},
        & $ 1 \leq i \leq n,$
        \label{eq:sum-inverse-lambda}
        \\
        \mu_i \coloneqq \frac1{a_{i+1} b_{i}-a_{i} b_{i+1} +c_{i+1}-c_{i}}\left(\frac{b_{i+1}-b_i}{v_{i+1}-v_{i}}\right)^2,
        & $ 1 \leq i \leq n - 1,$
        \label{eq:sum-inverse-mu}
                            \end{numcases}
    where  $ \mu_0 = \frac{a_1 x + c_1}{a_1^2 z^2},\, \mu_n \coloneqq 0$, and $b_0 \coloneqq x \neq b_1, v_0 \coloneqq z \neq 0$ are free parameters such that $a_1 x + c_1 \neq 0$.
\end{corollary}
\begin{proof}[Proof of Corollary \ref{corol:sum-inv}]
    Equation \eqref{eq:sum-inverse-lambda} is a restatement of equation \eqref{eq:single-pair-inverse-lambda}.  Equation \eqref{eq:sum-inverse-mu} follows from substituting equation \eqref{eq:th1-coefs} into equation \eqref{eq:single-pair-inverse-mu} $ \mu_i = \frac1{(v_{i+1}-v_{i})^2\beta_{i+1}} $.  Finally $\mu_0 \coloneqq \frac1{u_1 v_0 - u_0 v_1} $ with $u_0 \coloneqq 0, v_0 \coloneqq z $, and $ u_1 \coloneqq z \frac{a_1^2 z^2}{a_1 x + c_1} $ from equation \eqref{eq:th3-coefs-init}.
\end{proof}

Alternatively,
\begin{theorem}     \label{th:LSU-alt}
    Subject to the generalized continuant sequence $(\dot v_i)$ defined below satisfying $\dot v_i \neq \beta_i \dot v_{i-1}$, the sum of two single-pair matrices $\vec A \coloneqq \left(a_{\min(i,j)} b_{\max(i,j)} \right), \vec C \coloneqq \left(c_{\min(i,j)} \right)$ with $b_i\neq b_{i-1}$, admits the factorization
    \begin{gather*}
        \vec{A+C} = (\vec L\dot\Delta\vec B)\, \dot{\vec M}\, (\vec L\dot\Delta\vec B)^T,
        \qquad \dot \Delta \coloneqq \operatorname{diag} \left( \frac{b_i - b_{i-1}}{\dot v_i-\beta_i \dot v_{i-1}} \right),
        \qquad
        \dot{\vec M} \coloneqq \left(\dot u_{\min(i,j)} \dot v_{\max(i,j)} \right),
        \\
        \vec B \coloneqq  \begin{pmatrix}
        1 & & &  (0) \\
        - \beta_2 & 1  \\
        & \ddots & \ddots  \\
        (0) &  & -\beta_n & 1
        \end{pmatrix},
        \qquad \beta_i \coloneqq \frac{a_ib_{i-1}-a_{i-1}b_i + c_i - c_{i-1}}{(b_i-b_{i-1})^2}, \; 1\leq i \leq n,
    \end{gather*}
    with the generating pair $(\dot u_i), (\dot v_i)$ determined by the second-order coupled recursion
    \begin{equation}    \label{eq:th4-coefs-init}
        \left\{\begin{aligned}
            \dot u_1 & \coloneqq z\frac{a_1^2}{(b_1 - x)^2}
            \\
            \dot v_1 & \coloneqq z\frac{a_1b_1 + c_1}{(b_1 - x)^2} 
        \end{aligned}\right. ,        
    \end{equation}
    and, for $2 \leq i \leq n $,
    \begin{numcases}{}
            \dot u_i \coloneqq
             \frac{a_i b_{i-1}-a_{i-1} b_i +c_i-c_{i-1}}{(b_i - b_{i-1})^2} \dot u_{i-1} - 
                \dot u_{i-1} \frac{a_i - a_{i-1}}{b_i-b_{i-1}}
                + \dot v_{i-1} \left( \frac{a_i - a_{i-1}}{b_i-b_{i-1}} \right)^2
                \label{eq:th4-coefs-u}
            \\
            \dot v_i \coloneqq \left( 
                \frac{a_i b_i -2 a_{i-1} b_i +  a_{i-1} b_{i-1} +c_i -c_{i-1}}{\left(b_i-b_{i-1}\right)^2}-\frac{a_{i-1} b_{i-1}-2 a_{i-1} b_{i-2}+a_{i-2} b_{i-2}-(c_{i-1}-c_{i-2})}{\left(b_{i-1}-b_{i-2}\right)^2}
                \right) \dot v_{i-1}
            \nonumber\\ 
            \qquad -\, \frac{\left( a_{i-1} b_{i-2}-a_{i-2} b_{i-1} +c_{i-1}-c_{i-2} \right)^2 }{(b_{i-1}-b_{i-2})^4} \dot v_{i-2},
            \label{eq:th4-coefs-v}
    \end{numcases}
    where $ a_0 = c_0 = 0$ and $b_0 \coloneqq x \neq b_1,\; \dot v_0 \coloneqq z \neq 0$ are free parameters.
\end{theorem}
\begin{remark}
    $\vec L\dot{\Delta}\vec B$ is the lower triangular matrix
    \[
        \vec L\dot{\Delta}\vec B = \begin{pmatrix}
        \delta_1 \\
        \delta_1 - \beta_2\delta_2  & \delta_2 & & (0)  \\
        \delta_1 - \beta_2\delta_2  & \delta_2 - \beta_3 \delta_3 & \delta_3  \\
        \vdots & \vdots & \ddots & \ddots \\
        \delta_1 - \beta_2 \delta_2 & \delta_2 - \beta_3\delta_3 & \cdots & \delta_{n-1} - \beta_n \delta_n  & \delta_n
        \end{pmatrix},
        \qquad \delta_i \coloneqq \frac{b_i - b_{i-1}}{\dot v_i-\beta_i \dot v_{i-1}}.
    \]
\end{remark}
\begin{proof}
    From Theorem \ref{th:single-pair-sum}, $ \vec{A+C} = (\vec{LD_1 L})\, \vec T\, (\vec{LD_1 L})^T, 
    \ \vec D_1  \coloneqq \operatorname{diag}(b_1 - x, b_2 - b_1, \cdots, b_n - b_{n-1}) $.  From Corollary \ref{corol:tridiagonal-single-pair-2},
    $
        \vec T = \big( \vec L^{-1} \vec{\dot D}_2^{-1} \vec B \big)\: \vec{\dot M}\: \big( \vec L^{-1} \vec{\dot D}_2^{-1} \vec B \big)^T,
        \ \vec{\dot D}_2 \coloneqq \operatorname{diag}( \dot u_i - \beta_i \dot u_{i-1} )
    $.   Combining and simplifying yields the desired factorization, with
    \[
        \begin{dcases}
            \dot u_1 \coloneqq z (\alpha_1 - \beta _1 )^2
            \\
            \dot v_1 \coloneqq z \alpha_1
        \end{dcases},
    \]    
    and, for $2 \leq i \leq n $,
    \begin{numcases}{}
        \dot u_i \coloneqq \ \beta_i u_{i-1}  -  \dot u_{i-1} \left(\sum_{j=1}^i (\alpha_j - \beta_j - \beta_{j-1}) \right)
            + \dot v_{i-1} \left(\sum_{j=1}^i (\alpha_j - \beta_j - \beta_{j-1}) \right)^2
        \label{eq:th4-proof-recursion-u}\\
        \dot v_i \coloneqq \ \alpha_i \dot v_{i-1} -  \beta_{i-1}^2 \dot v_{i-2}.
        \label{eq:th4-proof-recursion-v}
    \end{numcases}
    Substituting equation \eqref{eq:th1-coefs} for $\alpha_i, \beta_i, \beta_{i-1}$ into the above expressions, we recover equations \eqref{eq:th4-coefs-init}, \eqref{eq:th4-coefs-u} and \eqref{eq:th4-coefs-v} after some algebra. 
\end{proof}

\begin{corollary}   \label{corol:sum-inv-alt}
    If $(\dot v_i)$ satisfies $\dot v_i \notin\{0, \beta_i \dot v_{i-1}\}$, then $\vec{A+C}$ is invertible in semi-closed form as
    \[
        (\vec{A+C})^{-1} = \big(\vec{B}^{-1}\dot\Delta^{-1}\vec L^{-1}\big)^T\, \dot{\vec M}^{-1}\, \big(\vec{B}^{-1}\dot\Delta^{-1}\vec L^{-1}\big),
        \qquad
        \dot{\vec M}^{-1} = \begin{psmallmatrix}
        \lambda_1 & -\mu_1 & & & (0) \\
        -\mu_1  & \lambda_2 & -\mu_2 \\
        & \ddots & \ddots & \ddots & \\
        & & -\mu_{n-2} & \lambda_{n-1} & -\mu_{n-1} \\
        (0) & & & -\mu_{n-1} & \lambda_n
        \end{psmallmatrix},
    \]
    with tridiagonal coefficients
    \begin{numcases}{}
        \lambda_i \coloneqq
         \frac{\dot v_{i+1}}{\dot v_i}\mu_{i} + \frac{\dot v_{i-1}}{\dot v_i}\mu_{i-1},
        & $ 1 \leq i \leq n,$
        \label{eq:sum-inverse-alt-lambda}
        \\
        \mu_i = \frac 1{\left(\dot v_{i+1}- \beta_{i+1} \dot v_{i}\right)^2},
        & $ 1 \leq i \leq n - 1,$
        \label{eq:sum-inverse-alt-mu}
                            \end{numcases}
    where $ \mu_0 \coloneqq \left(\frac{b_1 - x}{a_1 z}\right)^2,\, \mu_n \coloneqq 0$, and $b_0 \coloneqq x \neq b_1, \; \dot v_0 \coloneqq z \neq 0$ are free parameters.
\end{corollary}
\begin{remark}
    $\vec B^{-1}\dot{\Delta}^{-1}\vec L^{-1}$ is the lower triangular matrix
    \begin{align*}
        \vec B^{-1}\dot{\Delta}^{-1}\vec L^{-1}
        & = \begin{pmatrix}
        \frac 1{\delta_1} \\
        \frac {\beta_2}{\delta_1} - \frac 1{\delta_2}  & \frac 1{\delta_2} & & (0)  \\
        \frac {\beta_2\beta_3}{\delta_1} - \frac {\beta_3}{\delta_2}  & \frac {\beta_3}{\delta_2} - \frac 1{\delta_3} & \frac 1{\delta_3} \\
        \vdots & \vdots & \ddots & \ddots \\
        \frac {\beta_2\cdots\beta_n}{\delta_1} - \frac {\beta_3\cdots\beta_n}{\delta_2} & \frac {\beta_3\cdots\beta_n}{\delta_2} - \frac {\beta_4\cdots\beta_n}{\delta_3}
        & \cdots & \frac {\beta_n}{\delta_{n-1}} - \frac 1{\delta_n}  & \frac 1{\delta_n}
        \end{pmatrix}
        \\
        & =
        \left( \left(\frac{\prod_{k=j+1}^i \beta_k}{\delta_j} - \frac{\prod_{k=j+2}^i \beta_k}{\delta_{j+1}} \right)[ j < i ] + \frac1{\delta_i}[ i = j]
        \right)_{1\leq i,j\leq n},
    \end{align*}
    where any void product is deemed to be 1, $[\cdot]$ denotes Iverson's bracket, and $ \frac 1{\delta_i} \coloneqq \frac{\dot v_i-\beta_i \dot v_{i-1}}{b_i - b_{i-1}} $.
\end{remark}
\begin{proof}[Proof of Corollary \ref{corol:sum-inv-alt}]
    The factorization $(\vec{A+C})^{-1} = \big(\vec{B}^{-1}\dot\Delta^{-1}\vec L^{-1}\big)^T\, \dot{\vec M}^{-1}\, \big(\vec{B}^{-1}\dot\Delta^{-1}\vec L^{-1}\big)$ is a straighforward inversion of $\vec{A+C} = (\vec L\dot\Delta\vec B)\, \dot{\vec M}\, (\vec L\dot\Delta\vec B)^T$.  The calculation of $\dot{\vec M}^{-1}$ is given by Corollary \ref{corol:tridiagonal-single-pair-2}, equations \eqref{eq:single-pair-inverse-mu} and \eqref{eq:single-pair-inverse-lambda}, as
    \[
        \begin{dcases}
            \lambda_i \coloneqq
             \frac{\dot v_{i+1}}{\dot v_i}\mu_{i} + \frac{\dot v_{i-1}}{\dot v_i}\mu_{i-1},
            & 1 \leq i \leq n,
            \\
            \mu_i \coloneqq \frac1{\dot u_{i+1}\dot v_{i} - \dot u_{i} \dot v_{i+1}},
            &  1 \leq i \leq n - 1,
        \end{dcases}
    \]
    and the equation above for $\lambda_i$ matches \eqref{eq:sum-inverse-alt-lambda}.
    Substituting $ \dot u_i \coloneqq u_i\, \beta_1 \cdots  \beta_i,\ \dot v_i \coloneqq v_i\, \beta_1 \cdots \beta_i$ into the equation above for $\mu_i$ and simplifying,
    \[
        \mu_i = \frac 1{(\beta_1 \cdots \beta_i)^2 \beta_{i+1}} \frac1{u_{i+1}v_i - u_i v_{i+1}}
        = \frac 1{(\beta_1 \cdots \beta_{i+1})^2} \frac1{(v_{i+1} - v_i)^2}
    \]
    where the pair $(u_i), (v_i)$ is defined by equations \eqref{eq:th3-coefs-u} and \eqref{eq:th3-coefs-v}, and we plugged equation \eqref{eq:sum-inverse-mu} in the last step.  Substituting $ v_i = \dot v_i / \beta_1 \cdots \beta_i $ and simplifying, we recover equation \eqref{eq:sum-inverse-alt-mu} as required.  Finally, we have $\mu_0 \coloneqq \frac1{\dot u_1 \dot v_0 - \dot u_0 \dot v_1}$ with $\dot u_0 \coloneqq 0, \dot v_0 \coloneqq z$, and $\dot u_1 \coloneqq z\frac{a_1^2}{(b_1 - x)^2}$ from equation \eqref{eq:th4-coefs-init}.
\end{proof}

\begin{theorem}     \label{th5:SP-sum-inv}
    Subject to the generalized continuant sequence $(\dot v_i)$ defined in equations \eqref{eq:th4-coefs-init} and \eqref{eq:th4-coefs-v} satisfying $\dot v_i \notin\{0, \beta_i \dot v_{i-1}\}$, the sum of two single-pair matrices $\vec A \coloneqq \left(a_{\min(i,j)} b_{\max(i,j)} \right), \vec C \coloneqq \left(c_{\min(i,j)} \right)$ with $b_i\neq b_{i-1}$ is invertible with inverse coefficients
    \begin{flalign}
        (\vec {A+C})^{-1}_{j,j}
        =&\ \mathrlap{
        \r{[j<n-1]}\left(\frac{\dot v_{j+1}-\beta_{j+1} \dot v_{j}}{b_{j+1}-b_j} - \beta_{j+1}\frac{\dot v_j-\beta_j \dot v_{j-1}}{b_j-b_{j-1}}\right)^2 \sum_{k=j+2}^n \frac{ \beta_{j+2}^2\cdots\beta_{k-1}^2}{\dot v_k\dot v_{k-1}}
        + \frac1{(b_j-b_{j-1})^2}\frac{\dot v_{j-1}}{\dot v_j}
        } \nonumber \\
        &\ \mathrlap{
            +
            [j < n] \left(
                \frac{\dot v_j - \beta_j \dot v_{j-1}}{b_{j}-b_{j-1}}
                +
                \frac{\dot v_{j}}{b_{j+1} - b_{j}}
            \right)^2
            \frac{1}{\dot v_{j+1}\dot v_j},
        }
        & 1 \leq j \leq n;   \label{eq:suminv-diagonal}
    \end{flalign}\begin{flalign}
        (\vec {A+C})^{-1}_{i,j}
        =&\ \mathrlap{
        \r{[j<n-1]}\left(\frac{\dot v_{i+1}-\beta_{i+1} \dot v_{i}}{b_{i+1}-b_i} - \beta_{i+1}\frac{\dot v_i-\beta_i \dot v_{i-1}}{b_i-b_{i-1}}\right)\left(\frac{\dot v_{j+1}-\beta_{j+1} \dot v_{j}}{b_{j+1}-b_j} - \beta_{j+1}\frac{\dot v_j-\beta_j \dot v_{j-1}}{b_j-b_{j-1}}\right) 
        } \nonumber \\
        && \times \sum_{k=j+2}^n \frac{ \beta_{i+2}\cdots\beta_{k-1} \beta_{j+2}\cdots\beta_{k-1}}{\dot v_k\dot v_{k-1}}
        \nonumber
        \\
        &\ \mathrlap{
            + \frac{[i < j-1]}{b_j-b_{j-1}}\left(\frac{\dot v_{i+1}-\beta_{i+1} \dot v_{i}}{b_{i+1}-b_i} - \beta_{i+1}\frac{\dot v_i-\beta_i \dot v_{i-1}}{b_i-b_{i-1}}\right) \frac{\beta_{i+2}\cdots\beta_{j-1}}{\dot v_j}
        } \nonumber
        \\
        &\ \mathrlap{
            -[j<n]\left(\frac{\dot v_{i+1}-\beta_{i+1} \dot v_{i}}{b_{i+1}-b_i} - \beta_{i+1}\frac{\dot v_i-\beta_i \dot v_{i-1}}{b_i-b_{i-1}}\right)
            \left(
                \frac{\dot v_{j} - \beta_{j}\dot v_{j-1}}{b_{j}-b_{j-1}}
                + \frac{\dot v_{j}}{b_{j+1} - b_{j}}
            \right)
            \frac{\beta_{i+2}\cdots\beta_{j}}{\dot v_{j+1}\dot v_j}
        }   \nonumber
        \\
        &\ \mathrlap{
            -\frac{[i = j-1]}{b_j-b_{j-1}}\left(
                \frac{\dot v_{j-1} - \beta_{j-1}\dot v_{j-2}}{b_{j-1}-b_{j-2}}
                + \frac{\dot v_{j-1}}{b_j - b_{j-1}}
            \right) \frac{1}{\dot v_j},
        } 
        & 1 \leq i < j \leq n;   \label{eq:suminv-upper}
        \\
        (\vec {A+C})^{-1}_{i,j}
        =&\ \mathrlap{
        (\vec {A+C})^{-1}_{j,i}, \qquad\qquad 1 \leq j < i \leq n,
        }&
        \nonumber
    \end{flalign}    
    where $ a_0 = c_0 = 0, \beta_i \coloneqq \frac{a_ib_{i-1}-a_{i-1}b_i + c_i - c_{i-1}}{(b_i-b_{i-1})^2} $, and $b_0 \neq b_1 , \; \dot v_0 \coloneqq z \neq 0$ are free parameters.
\end{theorem}
\begin{proof}
    By Lemma \ref{lem:triang-prod}, %equation \eqref{eq:UTL-product},
    \[
        (\vec {A+C})^{-1}_{i,j}
        = \sum_{k=\max(i,j)}^n l_{k,i} \lambda_{k} l_{k,j} - \sum_{k=\max(i+1,j)}^n l_{k-1,i} \mu_{k-1} l_{k,j} - \sum_{k=\max(i,j+1)}^n l_{k,i} \mu_{k-1} l_{k-1,j}.
    \]
    Substituting $\lambda_k \coloneqq \frac{\dot v_{k+1}}{\dot v_k}\mu_k + \frac{\dot v_{k-1}}{\dot v_k}\mu_{k-1}$ with the convention $\mu_0 \coloneqq \frac{\left(b_1-x\right)^2}{(a_1 z)^2}, \mu_n \coloneqq 0$,
    \[
        (\vec {A+C})^{-1}_{i,j}
        =
        \sum_{k=\max(i,j)}^n l_{k,i} \left( \frac{\dot v_{k+1}}{\dot v_k}\mu_k + \frac{\dot v_{k-1}}{\dot v_k}\mu_{k-1} \right) l_{k,j} - \sum_{k=\max(i+1,j)}^n l_{k-1,i} \mu_{k-1} l_{k,j} - \sum_{k=\max(i,j+1)}^n l_{k,i} \mu_{k-1} l_{k-1,j}.
    \]    
    \textbullet\quad
    For off-diagonal coefficients where $1\leq i< j\leq n$, we have $\max(i,j) = \max(i+1,j) = j$ and $ \max(i,j+1) = j+1$. Splitting the first sum into two, reindexing one sum and gathering terms; factoring; and isolating the term for $k=j+1$,
    \begin{flalign}
        (\vec {A+C})^{-1}_{i,j}
        = &
        \mathrlap{\sum_{k=j+1}^n \left[ \left( \frac{\dot v_{k}}{\dot v_{k-1}}l_{k-1,j} - l_{k,j}\right) l_{k-1,i} 
        + \left( \frac{\dot v_{k-1}}{\dot v_k}l_{k,j} - l_{k-1,j}\right) l_{k,i} \right]  \mu_{k-1}}
        & \nonumber \\
        && + l_{j,i} \frac{\dot v_{j-1}}{\dot v_j} \mu_{j-1} l_{j,j} - l_{j-1,i}\mu_{j-1}l_{j,j},
        \nonumber \\
        = &
        -\mathrlap{\sum_{k=j+1}^n \left( \frac{\dot v_{k-1}}{\dot v_k}l_{k,j} - l_{k-1,j}\right) \left( \frac{\dot v_{k}}{\dot v_{k-1}}l_{k-1,i} - l_{k,i}\right) \mu_{k-1}
        + \left(\frac{\dot v_{j-1}}{\dot v_j}l_{j,i} - l_{j-1,i} \right) \mu_{j-1} l_{j,j}}
        & \nonumber \\
        = &
        -\mathrlap{\sum_{k=j+2}^n \underbrace{\left( \frac{\dot v_{k-1}}{\dot v_k}l_{k,j} - l_{k-1,j}\right) \left( \frac{\dot v_{k}}{\dot v_{k-1}}l_{k-1,i} - l_{k,i}\right) \mu_{k-1}}_{S_{i,j,k}}
        + \underbrace{\left(\frac{\dot v_{j-1}}{\dot v_j}l_{j,i} - l_{j-1,i} \right) \mu_{j-1} l_{j,j}}_{T_{i,j}}}
        & \nonumber \\ &&
        - \underbrace{\left( \frac{\dot v_{j}}{\dot v_{j+1}}l_{j+1,j} - l_{j,j}\right) \left( \frac{\dot v_{j+1}}{\dot v_{j}}l_{j,i} - l_{j+1,i}\right) \mu_{j}}_{U_{i,j}} [ j < n],
        \qquad\qquad 1\leq i< j\leq n.
        \label{eq:suminv-upper2}
    \end{flalign}
    Plugging $l_{i,j} \coloneqq \left(\frac{\beta_{j+1}}{\delta_j} - \frac{1}{\delta_{j+1}}\right) \beta_{j+2}\cdots\beta_i[j < i] + \frac1{\delta_i}[i = j]$ with appropriate index substitutions, we obtain after ad hoc factorizations
    \begin{flalign*}
        S_{i,j,k}
        & = \mathrlap{
        \left( \frac{\beta_{i+1}}{\delta_i} - \frac{1}{\delta_{i+1}} \right)\!\!\left( \frac{\beta_{j+1}}{\delta_j} - \frac{1}{\delta_{j+1}} \right)\!\! \left( \frac{\dot v_{k-1}}{\dot v_k}\beta_k - 1\right)\!\! \left( \frac{\dot v_{k}}{\dot v_{k-1}} - \beta_k\right) \mu_{k-1} \beta_{i+2}\cdots\beta_{k-1} \beta_{j+2}\cdots\beta_{k-1},
        }&
        \\ &&
        1 \leq i < j < k - 1 \leq n - 1;
        \\
        T_{i,j}
        & = \mathrlap{
        \begin{dcases}
            \left(\frac{\beta_{i+1}}{\delta_i} - \frac{1}{\delta_{i+1}}\right) \left(\frac{\dot v_{j-1}}{\dot v_j}\beta_j - 1 \right) \frac{\mu_{j-1}}{\delta_j} \beta_{i+2}\cdots\beta_{j-1},
            & 1\leq i< j-1\leq n-1,
            \\
            \left(
                \frac{\dot v_{j-1}}{\dot v_j}\left(\frac{\beta_j}{\delta_{j-1}} - \frac1{\delta_j} \right)-\frac1{\delta_{j-1}}
            \right) \frac{\mu_{j-1}}{\delta_j},
            & 1\leq i = j-1 \leq n - 1;
        \end{dcases}
        } &
        \\
        U_{i,j}
        & = \mathrlap{
            \left(\frac{\beta_{i+1}}{\delta_i} - \frac{1}{\delta_{i+1}}\right)
            \left(
                \frac{\dot v_{j+1}}{\dot v_{j}} - \beta_{j+1}
            \right)
            \left(
                \frac{\dot v_{j}}{\dot v_{j+1}}\left(\frac{\beta_{j+1}}{\delta_{j}} - \frac1{\delta_{j+1}} \right)-\frac1{\delta_{j}}
            \right) \mu_{j} \beta_{i+2}\cdots\beta_{j},
        }&
        \\ &&
        1\leq i < j \leq n - 1.
    \end{flalign*}
    Plugging $\mu_j \coloneqq 1/(\dot v_{j+1} - \beta_{j+1}\dot v_{j})^2$ and similarly for $\mu_{k-1}, \mu_{j-1}$, and simplifying,
    \begin{flalign*}
        S_{i,j,k}
        & = \mathrlap{
        -\left( \frac{\beta_{i+1}}{\delta_i} - \frac{1}{\delta_{i+1}} \right)\!\!\left( \frac{\beta_{j+1}}{\delta_j} - \frac{1}{\delta_{j+1}} \right)\frac{ \beta_{i+2}\cdots\beta_{k-1} \beta_{j+2}\cdots\beta_{k-1}}{\dot v_k\dot v_{k-1}},}
        & 1 \leq i < j < k - 1 \leq n - 1;
        \\
        T_{i,j}
        & = \mathrlap{
        \begin{dcases}
            -\left(\frac{\beta_{i+1}}{\delta_i} - \frac{1}{\delta_{i+1}}\right) \frac{\beta_{i+2}\cdots\beta_{j-1}}{\delta_j \dot v_j (\dot v_j - \beta_j\dot v_{j-1})},
            & 1\leq i< j-1\leq n-1,
            \\
            -\left(
                \frac1{\delta_{j-1}}
                + \frac{\dot v_{j-1}}{\delta_{j}(\dot v_{j} - \beta_{j}\dot v_{j-1})}
            \right) \frac{1}{\delta_j\dot v_j (\dot v_{j} - \beta_{j}\dot v_{j-1})},
            & 1\leq i = j-1 \leq n - 1;
        \end{dcases}
        } &
        \\
        U_{i,j}
        & = \mathrlap{
            -\left(\frac{\beta_{i+1}}{\delta_i} - \frac{1}{\delta_{i+1}}\right)
            \left(
                \frac1{\delta_{j}}
                + \frac{\dot v_{j}}{\delta_{j+1}(\dot v_{j+1} - \beta_{j+1}\dot v_{j})}
            \right) \frac{\beta_{i+2}\cdots\beta_{j}}{\dot v_{j+1}\dot v_j},
        }
        & 1\leq i < j \leq n - 1.
    \end{flalign*}    
    Finally, plugging $ \frac 1{\delta_i} \coloneqq \frac{\dot v_i-\beta_i \dot v_{i-1}}{b_i - b_{i-1}} $ with ad hoc index substitutions,
    \begin{flalign*}
        S_{i,j,k}
        & = \mathrlap{
        -\left(\frac{\dot v_{i+1}-\beta_{i+1} \dot v_{i}}{b_{i+1}-b_i} - \beta_{i+1}\frac{\dot v_i-\beta_i \dot v_{i-1}}{b_i-b_{i-1}}\right)\left(\frac{\dot v_{j+1}-\beta_{j+1} \dot v_{j}}{b_{j+1}-b_j} - \beta_{j+1}\frac{\dot v_j-\beta_j \dot v_{j-1}}{b_j-b_{j-1}}\right)}&
        \\ &&
        \times\frac{ \beta_{i+2}\cdots\beta_{k-1} \beta_{j+2}\cdots\beta_{k-1}}{\dot v_k\dot v_{k-1}},
        \qquad\qquad
        1 \leq i < j < k - 1 \leq n - 1;
        \\
        T_{i,j}
        & = \mathrlap{
        \begin{dcases}
            \frac1{b_j-b_{j-1}}\left(\frac{\dot v_{i+1}-\beta_{i+1} \dot v_{i}}{b_{i+1}-b_i} - \beta_{i+1}\frac{\dot v_i-\beta_i \dot v_{i-1}}{b_i-b_{i-1}}\right) \frac{\beta_{i+2}\cdots\beta_{j-1}}{\dot v_j},
            & 1\leq i< j-1\leq n-1,
            \\
            -\frac1{b_j-b_{j-1}}\left(
                \frac{\dot v_{j-1} - \beta_{j-1}\dot v_{j-2}}{b_{j-1}-b_{j-2}}
                + \frac{\dot v_{j-1}}{b_j - b_{j-1}}
            \right) \frac{1}{\dot v_j},
            & 1\leq i = j-1 \leq n - 1;
        \end{dcases}
        } &
        \\
        U_{i,j}
        & = \mathrlap{
            \left(\frac{\dot v_{i+1}-\beta_{i+1} \dot v_{i}}{b_{i+1}-b_i} - \beta_{i+1}\frac{\dot v_i-\beta_i \dot v_{i-1}}{b_i-b_{i-1}}\right)
            \left(
                \frac{\dot v_{j} - \beta_{j}\dot v_{j-1}}{b_{j}-b_{j-1}}
                + \frac{\dot v_{j}}{b_{j+1} - b_{j}}
            \right)
            \frac{\beta_{i+2}\cdots\beta_{j}}{\dot v_{j+1}\dot v_j},
        }&
        \\
        && 1\leq i < j \leq n - 1.
    \end{flalign*}    
    Substituting the above expressions into equation \eqref{eq:suminv-upper2}, we recover equation \eqref{eq:suminv-upper} as required. 
    \\
    \textbullet\quad 
    For diagonal coefficients where $i=j$, we have
    \[
        (\vec {A+C})^{-1}_{j,j}
        =
        \sum_{k=j}^n \left( \frac{\dot v_{k+1}}{\dot v_k}\mu_k + \frac{\dot v_{k-1}}{\dot v_k}\mu_{k-1} \right) l_{k,j}^2 - 2\sum_{k=j+1}^n l_{k,j} \mu_{k-1} l_{k-1,j}.
    \]
    Splitting the sum, reindexing as appropriate and gathering terms,
    \begin{flalign}
        (\vec {A+C})^{-1}_{j,j}
        & = \mathrlap{
        \sum_{k=j+1}^n \left( \frac{\dot v_{k}}{\dot v_{k-1}}l_{k-1,j}^2 + \frac{\dot v_{k-1}}{\dot v_k} l_{k,j}^2 - 2 l_{k,j}l_{k-1,j} \right)\mu_{k-1} + \frac{\dot v_{j-1}}{\dot v_j} l_{j,j}^2\mu_{j-1}
        }& \nonumber \\
        & = \mathrlap{
        \sum_{k=j+1}^n \left( l_{k-1,j} - \frac{\dot v_{k-1}}{\dot v_k} l_{k,j}\right)^2 \frac{\dot v_{k}}{\dot v_{k-1}}\mu_{k-1} + \frac{\dot v_{j-1}}{\dot v_j} l_{j,j}^2\mu_{j-1}
        }& \nonumber \\
        & = \mathrlap{
        \sum_{k=j+2}^n \underbrace{\left( l_{k-1,j} - \frac{\dot v_{k-1}}{\dot v_k} l_{k,j}\right)^2 \frac{\dot v_{k}}{\dot v_{k-1}}\mu_{k-1}}_{S_{j,j,k}}
        + \underbrace{\frac{\dot v_{j-1}}{\dot v_j} l_{j,j}^2\mu_{j-1}}_{T_{j,j}}
        }& \nonumber \\
        && + \underbrace{\left( l_{j,j} - \frac{\dot v_{j}}{\dot v_{j+1}} l_{j+1,j}\right)^2 \frac{\dot v_{j+1}}{\dot v_{j}}\mu_j}_{U_{j,j}} [j < n],
        \qquad 1 \leq j \leq n.
        \label{eq:suminv-diagonal2}
    \end{flalign}
    Following similar steps as before, we find
    \begin{align*}
        S_{j,j,k}
        & =         \left(\frac{\dot v_{j+1}-\beta_{j+1} \dot v_{j}}{b_{j+1}-b_j} - \beta_{j+1}\frac{\dot v_j-\beta_j \dot v_{j-1}}{b_j-b_{j-1}}\right)^2 \frac{ \beta_{j+2}^2\cdots\beta_{k-1}^2}{\dot v_k\dot v_{k-1}},
                & 1 \leq j < k - 1 \leq n - 1;
        \\
        T_{j,j}
        & =             \frac1{(b_j-b_{j-1})^2}\frac{\dot v_{j-1}}{\dot v_j},
                & 1\leq j \leq n;
        \\
        U_{j,j}
        & =             \left(
                \frac{\dot v_j - \beta_j \dot v_{j-1}}{b_{j}-b_{j-1}}
                +
                \frac{\dot v_{j}}{b_{j+1} - b_{j}}
            \right)^2
            \frac{1}{\dot v_{j+1}\dot v_j},
                & 1 \leq j \leq n - 1.
    \end{align*}
    Substituting the above expressions into equation \eqref{eq:suminv-diagonal2}, we recover equation \eqref{eq:suminv-diagonal} as required. 
\end{proof}

\section{Application: symbolic inverse of the Gram matrix of a system of ramp functions}        \label{sec:Gram-inv}

In the Hilbert space of continuous functions over the interval $[0, 1]$, let $f_i(x) \coloneqq \max(0, k_i - x)$ be a system of $n$ ramp functions with shift parameters $0 < k_1 < \cdots < k_n \leq 1$.  The Gram matrix $\vec G$ of the system has coefficients
\[
    g_{i,j} \coloneqq \inner{f_i,f_j}
    \coloneqq
        \int_0^1 \max(0, k_i - x) \max(0, k_j - x) \d x
    =
        \frac12 k_{\min(i,j)}^2 k_{\max(i,j)} - \frac 16 k_{\min(i,j)}^3,
\]
so that $\vec G = \frac16(\vec A + \vec C)$ where $\vec A + \vec C$ is the sum of two single-pair matrices with generators $a_i \coloneqq 3 k_i^2, b_i \coloneqq k_i, c_i \coloneqq -k_i^3$.  For $n=2$ we have
\[
    6\vec G
    =
    \begin{pmatrix}
        2k_1^3            &  k_1^2(3k_2 - k_1)
        \\
        k_1^2(3k_2 - k_1) &  2 k_2^3
    \end{pmatrix}
    , \quad
    \vec G^{-1}
    =
    \frac 6{k_1^3 (k_2 - k_1)^2 (4k_2-k_1)}
    \begin{pmatrix}
        2k_2^3             &  -k_1^2(3k_2 - k_1)
        \\
        -k_1^2(3k_2 - k_1) &  2 k_1^3
    \end{pmatrix}
\]

\subsection{Factorizations}

By Theorem \ref{th:single-pair-sum} for $x = 0$,
\begin{gather*}
    6\vec{G} = (\vec{LDL})\, \vec T\, (\vec{LDL})^T,
    \quad \vec D  \coloneqq \operatorname{diag}(k_1, k_2 - k_1, \cdots, k_n - k_{n-1}),
    \nonumber\\
    \vec T \coloneqq \begin{pmatrix}
    \alpha_1 & -\beta_1 & & & (0) \\
    -\beta_1  & \alpha_2 & -\beta_2 \\
    & \ddots & \ddots & \ddots & \\
    & & -\beta_{n-2} & \alpha_{n-1} & -\beta_{n-1} \\
    (0) & & & -\beta_{n-1} & \alpha_n
    \end{pmatrix},
\end{gather*}
with tridiagonal coefficients that remarkably simplify to
\begin{equation}    \label{eq:gram-th1}
    \begin{dcases}
        \alpha_1 = 2 k_1
        \\
        \beta_1 = -k_1
    \end{dcases},
    \quad
    \begin{dcases}
        \alpha_i \coloneqq \beta_{i-1}+\beta_i + \frac{3k_i^2 - 3k_{i-1}^2}{k_i-k_{i-1}} - \frac{3k_{i-1}^2 - 3k_{i-2}^2}{k_{i-1}-k_{i-2}} = 2(k_i-k_{i-2})
        \\
        \beta_i \coloneqq \frac{3 k_i^2 k_{i-1} - 3k_{i-1}^2k_i - (k_i^3 - k_{i-1}^3)}{(k_i - k_{i-1})^2}
        = -(k_i-k_{i-1})
    \end{dcases},
    \quad 2 \leq i \leq n
\end{equation}
where $k_0 \coloneqq 0$.  Theorem \ref{th:LSU} for $x = 0, z = 1$ gives the factorization
\[
    6\vec{G}= (\vec{L\Delta L}^{-1})\, \vec M\, (\vec{L\Delta L}^{-1})^T,
    \qquad \vec \Delta = \operatorname{diag} \left( \frac{k_i - k_{i-1}}{v_i-v_{i-1}} \right),
    \qquad \vec M \coloneqq \left(u_{\min(i,j)} v_{\max(i,j)} \right),
\]
with the generating pair $(u_i), (v_i)$ determined by the second-order coupled recursion
\begin{equation*}
    \left\{\begin{aligned}        u_1 & = \frac{a_1^2}{c_1} = -9k_1
        \\
        v_1 & = 1 + \frac{a_1b_1}{c_1} = -2
    \end{aligned}\right. ,        
\end{equation*}
and, for $2 \leq i \leq n $ (equations \eqref{eq:th3-proof-recursion-u} and \eqref{eq:th3-proof-recursion-v}),
\[ \begin{dcases}
    u_i \coloneqq \ u_{i-1}  -  \frac {u_{i-1}}{\beta_i} \left(\sum_{j=1}^i (\alpha_j - \beta_j - \beta_{j-1}) \right)
        + \frac{v_{i-1}}{\beta_i}\left(\sum_{j=1}^i (\alpha_j - \beta_j - \beta_{j-1}) \right)^2
    \\
    v_i \coloneqq \ \frac{\alpha_i}{\beta_i} v_{i-1} -  \frac{\beta_{i-1}}{\beta_i}{v_{i-2}}
    \end{dcases}
\]
where $\beta_0 \coloneqq 0, v_0 \coloneqq 1$.  Substituting equation \eqref{eq:gram-th1} and simplifying,
\[ \begin{dcases}
        u_i =\ u_{i-1}  +  3\frac {k_i + k_{i-1}}{k_i-k_{i-1}} u_{i-1}
        - 9\frac{\left(k_i + k_{i-1}\right)^2}{k_i-k_{i-1}}v_{i-1}
        \\
        v_i =\ -2\frac{k_i - k_{i-2}}{k_i - k_{i-1}} v_{i-1}
        - \frac{k_{i-1} - k_{i-2}}{k_i - k_{i-1}}{v_{i-2}}
    \end{dcases},
    \qquad 2 \leq i \leq n.
\]

Theorem \ref{th:LSU-alt} for $x = 0, z = 1$ gives the factorization 
\begin{gather*}
    6\vec{G} = (\vec L\dot\Delta\vec B)\, \dot{\vec M}\, (\vec L\dot\Delta\vec B)^T,
    \qquad \dot \Delta = \operatorname{diag} \left( \frac{k_i - k_{i-1}}{\dot v_i+(k_i-k_{i-1}) \dot v_{i-1}} \right),
    \qquad
    \dot{\vec M} \coloneqq \left(\dot u_{\min(i,j)} \dot v_{\max(i,j)} \right),
    \\
    \vec B =  \begin{pmatrix}
    1 & & &  (0) \\
    k_2 - k_1 & 1  \\
    & \ddots & \ddots  \\
    (0) &  & k_n - k_{n-1} & 1
    \end{pmatrix},
\end{gather*}
with the generating pair $(\dot u_i), (\dot v_i)$ determined by the second-order coupled recursion
\[        \left\{\begin{aligned}        \dot u_1 & = 9 k_1^2
        \\
        \dot v_1 & = 2 k_1
    \end{aligned}\right. ,        
\]
and, for $2 \leq i \leq n $,
\[
    \begin{dcases}
        \dot u_i =
         -(k_i - k_{i-1}) \dot u_{i-1} - 
            3\dot u_{i-1} (k_i+k_{i-1})
            + 9 \dot v_{i-1} (k_i+k_{i-1})^2
                    \\
        \dot v_i = 2(k_i - k_{i-2})\dot v_{i-1}
        - (k_{i-1} - k_{i-2})^2 \dot v_{i-2}
            \end{dcases}
\]
where $b_0 \coloneqq 0 ,\dot v_0 \coloneqq 1$.

\subsection{Symbolic inverses}

Corollary \ref{corol:sum-inv-alt} gives the inverse of $6\vec{G}$ as
\[
    (6\vec{G})^{-1} = \big(\vec{B}^{-1}\dot\Delta^{-1}\vec L^{-1}\big)^T\, \dot{\vec M}^{-1}\, \big(\vec{B}^{-1}\dot\Delta^{-1}\vec L^{-1}\big),
\]
where $\dot{\vec M}^{-1}$ is the symmetric tridiagonal matrix
\[
    \dot{\vec M}^{-1} = \begin{pmatrix}
    \lambda_1 & -\mu_1 & & & (0) \\
    -\mu_1  & \lambda_2 & -\mu_2 \\
    & \ddots & \ddots & \ddots & \\
    & & -\mu_{n-2} & \lambda_{n-1} & -\mu_{n-1} \\
    (0) & & & -\mu_{n-1} & \lambda_n
    \end{pmatrix},
\]
with tridiagonal coefficients
\[
    \begin{dcases}
        \lambda_i \coloneqq
         \frac{\dot v_{i+1}}{\dot v_i}\mu_{i} + \frac{\dot v_{i-1}}{\dot v_i}\mu_{i-1},
        & 1 \leq i \leq n,
        \\
        \mu_i = \frac 1{\left(\dot v_{i+1} + (k_{i+1} - k_i) \dot v_{i}\right)^2},
        & 1 \leq i \leq n - 1,
    \end{dcases}
\]
where $\mu_0 = \frac1{9 k_1^2}, \mu_n = b_0 = 0,\dot v_0 = 1$.  In addition, $\vec B^{-1}\dot{\Delta}^{-1}\vec L^{-1}$ is the lower triangular matrix
\[
    \left(\vec B^{-1}\dot{\Delta}^{-1}\vec L^{-1}\right)_{i,j}
    =
    (-1)^{i-j} \left(\frac{k_{j+1} - k_j}{\delta_j} + \frac{1}{\delta_{j+1}} \right)\left( \prod_{m=j+2}^i (k_m - k_{m-1}) \right)[ j < i ] + \frac1{\delta_i}[ i = j],
\]
where any void product is deemed to be 1, $[\cdot]$ denotes Iverson's bracket, and
\[
    \frac 1{\delta_i} = \frac{\dot v_i + (k_i - k_{i-1}) \dot v_{i-1}}{k_i - k_{i-1}} = \dot v_{i-1} + \frac{\dot v_i}{k_i - k_{i-1}}.
\]

For ease of writing, we introduce the notation $x_i \coloneqq b_i - b_{i-1} = k_i - k_{i-1}$.  By Theorem \ref{th5:SP-sum-inv},
\begin{flalign*}
    (6\vec{G})^{-1}_{j,j}
    =& \mathrlap{
    \left(\frac{\dot v_{j+1}-\beta_{j+1} \dot v_{j}}{x_{j+1}} - \beta_{j+1}\frac{\dot v_j-\beta_j \dot v_{j-1}}{x_j}\right)^2 
    \sum_{k=j+2}^n \frac{ \beta_{j+2}^2\cdots\beta_{k-1}^2}{\dot v_k\dot v_{k-1}}
    + \frac1{x_j^2}\frac{\dot v_{j-1}}{\dot v_j}
    }     \\
    &&    
        +
        \frac{[j < n]}{\dot v_{j+1}\dot v_j} \left(
            \frac{\dot v_j - \beta_j\dot v_{j-1}}{x_j} + \frac{\dot v_j}{x_{j+1}}
        \right)^2,
    \qquad 1 \leq j \leq n;
\end{flalign*}
\begin{flalign*}
    (6\vec G)^{-1}_{i,j}
    =&\ \mathrlap{
    \left(\frac{\dot v_{i+1}-\beta_{i+1} \dot v_{i}}{x_{i+1}} - \beta_{i+1}\frac{\dot v_i-\beta_i \dot v_{i-1}}{x_i}\right)\left(\frac{\dot v_{j+1}-\beta_{j+1} \dot v_{j}}{x_{j+1}} - \beta_{j+1}\frac{\dot v_j-\beta_j \dot v_{j-1}}{x_j}\right) 
    }     \\
    && \times \sum_{k=j+2}^n \frac{ \beta_{i+2}\cdots\beta_{k-1} \beta_{j+2}\cdots\beta_{k-1}}{\dot v_k\dot v_{k-1}}
        \\
    &\ \mathrlap{
        + \frac{[i < j-1]}{x_j}\left(\frac{\dot v_{i+1}-\beta_{i+1} \dot v_{i}}{x_{i+1}} - \beta_{i+1}\frac{\dot v_i-\beta_i \dot v_{i-1}}{x_i}\right) \frac{\beta_{i+2}\cdots\beta_{j-1}}{\dot v_j}
    }     \\
    &\ \mathrlap{
        -[j<n]\left(\frac{\dot v_{i+1}-\beta_{i+1} \dot v_{i}}{x_{i+1}} - \beta_{i+1}\frac{\dot v_i-\beta_i \dot v_{i-1}}{x_i}\right)
        \left(
            \frac{\dot v_{j} - \beta_{j}\dot v_{j-1}}{x_{j}}
            + \frac{\dot v_{j}}{x_{j+1}}
        \right)
        \frac{\beta_{i+2}\cdots\beta_{j}}{\dot v_{j+1}\dot v_j}
    }       \\
    &\ \mathrlap{
        -\frac{[i = j-1]}{x_j}\left(
            \frac{\dot v_{j-1} - \beta_{j-1}\dot v_{j-2}}{x_{j-1}}
            + \frac{\dot v_{j-1}}{x_j}
        \right) \frac{1}{\dot v_j},
    } 
    & 1 \leq i < j \leq n,   \nonumber\end{flalign*}    
where $\dot v_i$ is given by the recurrence equation  $ \dot v_i = 2(x_i + x_{i-1}) \dot v_{i-1} - x_{i-1}^2 \dot v_{i-2}$. Plugging $ \beta_i = -x_i $ with ad hoc index substitutions and simplifying,
\begin{flalign*}
    (6\vec G)^{-1}_{j,j}
    =&\ \mathrlap{
    \left(\frac{\dot v_{j+1}}{x_{j+1}} + \frac{x_{j+1}+x_j}{x_j}\dot v_{j} + x_{j+1}\dot v_{j-1} \right)^2  \sum_{k=j+2}^n \frac{ x_{j+2}^2\cdots x_{k-1}^2}{\dot v_k\dot v_{k-1}}
    + \frac1{x_j^2}\frac{\dot v_{j-1}}{\dot v_j}
    }     \\
    &\ \mathrlap{   
        + 
        \frac{[j < n]}{\dot v_{j+1}\dot v_j} \left(
            \left(\frac{1}{x_{j+1}}+\frac{1}{x_j}\right)\dot v_j+\dot v_{j-1}
        \right)^2,
    }
    & 1 \leq j \leq n;       \\
    (6\vec G)^{-1}_{i,j}
    =&\ \mathrlap{
    (-1)^{j-i} \left(\frac{\dot v_{i+1}}{x_{i+1}} + \frac{x_{i+1}+x_i}{x_i}\dot v_{i} + x_{i+1}\dot v_{i-1} \right) \left(\frac{\dot v_{j+1}}{x_{j+1}} + \frac{x_{j+1}+x_j}{x_j}\dot v_{j} + x_{j+1}\dot v_{j-1} \right)
    }     \\
    && \times \sum_{k=j+2}^n \frac{ x_{i+2}\cdots x_{k-1} x_{j+2}\cdots x_{k-1}}{\dot v_k\dot v_{k-1}}
        \\
    &\ \mathrlap{
        + [i < j-1] \frac{(-1)^{j-i}}{x_j}\left(\frac{\dot v_{i+1}}{x_{i+1}} + \frac{x_{i+1}+x_i}{x_i}\dot v_{i} + x_{i+1}\dot v_{i-1} \right) \frac{x_{i+2}\cdots x_{j-1}}{\dot v_j}
    }     \\
    &\ \mathrlap{
        +[j<n](-1)^{j-i}\left(\frac{\dot v_{i+1}}{x_{i+1}} + \frac{x_{i+1}+x_i}{x_i}\dot v_{i} + x_{i+1}\dot v_{i-1} \right)
        \left(
             \left(\frac{1}{x_{j+1}}+\frac{1}{x_j}\right)\dot v_j+\dot v_{j-1}
        \right)
        \frac{x_{i+2}\cdots x_{j}}{\dot v_{j+1}\dot v_j}
    }       \\
    &\ \mathrlap{
        -\frac{[i = j-1]}{x_j}\left(
            \left(\frac{1}{x_j} + \frac{1}{x_{j-1}}\right)\dot v_{j-1} + \dot v_{j-2}
        \right) \frac{1}{\dot v_j},
    } 
    & 1 \leq i < j \leq n.   \nonumber\end{flalign*}    
If so desired, expressions of the form $\frac{\dot v_{i+1}}{x_{i+1}} + \frac{x_{i+1}+x_i}{x_i}\dot v_{i} + x_{i+1}\dot v_{i-1}$ may be rewritten by means of the recurrence equation as
\[
    \frac{\dot v_{i+1}}{x_{i+1}} + \frac{x_{i+1}+x_i}{x_i}\dot v_{i} + x_{i+1}\dot v_{i-1}
    =
    \frac{x_{i+1}+x_i}{x_{i+1}}\left[\left(2 + \frac{x_{i+1}}{x_i}\right)\dot v_{i} + (x_{i+1}-x_i)\dot v_{i-1} \right],
\]
which removes dependence on $\dot v_{i+1}$.

\subsection{Continuant sequence}

It is easy to see that the continuants $(\dot v_i)$ form a sequence of homogeneous integer polynomials of incremental degree and number of variables:
\[
    \dot v_i = \sum_{\substack{(p_1, \cdots, p_i)\in\{0,1,2\}^i \\ p_1+\cdots+p_i = i}} \gamma({p_1, \dots, p_i})\, x_1^{p_1} \cdots x_i^{p_i},
    \qquad \gamma({p_1, \dots, p_i}) \in\mathbb Z.
\]
Substituting the above expression into the recurrence equation $ \dot v_{i+1} = 2x_{i+1} \dot v_i + x_i (2\dot v_i - x_i\dot v_{i-1}) $ and inspecting coefficients, we obtain the coefficient recurrence equations
\[
    \left\{\begin{array}{llll}
        \gamma({p_1, \dots, p_{i-1}, p_i, 1})
        & = &
        2\gamma({p_1, \dots, p_{i-1}, p_i}), & p_{i+1} = 1, p_i \in \{0,1\},
        \\[\medskipamount]
        \gamma({p_1, \dots, p_{i-1}, p_i + 1, 0})
        & = &
        2\gamma({p_1, \dots, p_{i-1}, p_i}) - \gamma(p_1, \dots, p_{i-1}), & p_{i+1} = 0, p_i \in \{0,1\},
        \\[\medskipamount]
        \gamma({p_1, \dots, p_{i-1}, p_i, p_{i+1}}) & = & 0, & \text{otherwise.}
    \end{array}\right.
\]
Further, we state without proof the following polynomial coefficient bounds:
\[\arraycolsep=1.4pt
    \gamma({p_1, \dots, p_i}) = 0 \quad\text{or}\quad
    \left\{\begin{array}{llll}
        3^{i/2} & \leq \gamma({p_1, \dots, p_i}) & \leq 2^i,      & \qquad i\text{ even},    
        \\[\medskipamount]
        2 \times 3^{(i-1)/2} & \leq \gamma({p_1, \dots, p_i}) & \leq 2^i,  & \qquad i\text{ odd}.
    \end{array}\right.
\]

Given their connection with continued fractions, the polynomials $(\dot v_i)$ are likely orthogonal with respect to some measure or moment functional \citep[see][for theoretical background]{dunkl_xu:2014}.  

\section{Application: numerical inverse of the sum of two single-pair matrices}

The semi-closed form formulas \eqref{eq:suminv-diagonal} and \eqref{eq:suminv-upper} based on the continuant sequence $(\dot v_i)$ given by equations \eqref{eq:th4-coefs-init} and \eqref{eq:th4-coefs-v} readily translate into a numerical computation algorithm for $(\vec A + \vec C)^{-1}$.  Algorithm 1 below provides the corresponding pseudocode, which is implemented in C in Appendix.  Inspection of the pseudocode reveals two sources of potential numerical instability due to division near zero:
\begin{itemize}[leftmargin=*,parsep=3pt,topsep=-2pt]
    \item When any two successive $b_i,b_{i+1}$ inputs are close.  We do not examine this case as it is inherent to inversion formulas \eqref{eq:suminv-diagonal} and \eqref{eq:suminv-upper}; however, this limitation could in principle be resolved mathematically by calculating limits of the form $\lim_{b_{i+1} \to b_i}$ (which exist since the left-hand side coefficient $(\vec{A+C})_{i,j}^{-1}$ is finite.)
    \item When any continuant value $\dot v_k$ is small, which is equivalent to say that the matrix $\vec{A+C}$ is ill-conditioned.
\end{itemize}
We empirically studied the numerical stability of Algorithm 1 against the QR method for ill-conditioned $3\times3$ matrices using two methodologies:
\begin{enumerate}[leftmargin=*]
    \item Using specific inputs $\vec a = (1, 1, 1)^T, \vec b = (1,5/3,3)^T, \vec c = (0,1,\varepsilon-3)^T, x = 0, z = 1$ for which $\det(\vec{A+C}) = -\frac{\varepsilon}9$, for $\varepsilon > 0$.  In this approach, the inverse $(\vec{A+C})^{-1}$ is available in closed form and was used to measure the accuracy of the algorithm output.
    \item Using systematic inputs $\vec a, \vec c$ together with fixed inputs $\vec b = (1,-1,1)^T, x = 0, z = 1$ so that $\operatorname{Sp}(\vec{A+C}) = \{-1, \varepsilon, 1\}$, for $\varepsilon = 10^{-1}, 10^{-4}, 10^{-7}$.  In this approach, we use Cramer's rule to measure for accuracy.
\end{enumerate}

\begin{algorithm}[!htbp]
\small
%\LinesNumbered
\DontPrintSemicolon
\SetAlgoLined
\SetKwInOut{AlgName}{name}
\SetKwInOut{Input}{input}\SetKwInOut{Output}{output}
%\SetKw{Let}{let}
\AlgName{SPSumInverse}
\Input{Size $n$, generators $\vec{a,b,c}\in\mathbb R^n$, free parameters $x,z$}
\Output{Upper-triangular matrix coefficients $(q_{i,j})_{1\leq i\leq j\leq n}$ } 
\;
Initialize $ (\dot v_i)_{0\leq i\leq n}, (\beta_i)_{1\leq i\leq n}, (\pi_{i,j})_{2\leq i\leq j+1\leq n}, (q_{i,j})_{1\leq i\leq j\leq n} $
\;
$ b_0 \coloneqq x,\ \dot v_0 \coloneqq z,\ \dot v_1 \coloneqq z\frac{a_1b_1 + c_1}{(b_1 - x)^2},\ \beta _1=\frac{a_1 x+c_1}{\left(b_1-x\right){}^2}$ 
\;
\For{i=2 \KwSty{to} n}{
    $ \begin{aligned}
        \dot v_i \coloneqq & \left( 
        \frac{a_i b_i -2 a_{i-1} b_i +  a_{i-1} b_{i-1} +c_i -c_{i-1}}{\left(b_i-b_{i-1}\right)^2}-\frac{a_{i-1} b_{i-1}-2 a_{i-1} b_{i-2}+a_{i-2} b_{i-2}-(c_{i-1}-c_{i-2})}{\left(b_{i-1}-b_{i-2}\right)^2}
        \right) \dot v_{i-1}
        \\ 
        & -\, \frac{\left( a_{i-1} b_{i-2}-a_{i-2} b_{i-1} +c_{i-1}-c_{i-2} \right)^2 }{(b_{i-1}-b_{i-2})^4} \dot v_{i-2}
        \\
        \beta_i \coloneqq &\ \frac{a_ib_{i-1}-a_{i-1}b_i + c_i - c_{i-1}}{(b_i-b_{i-1})^2}
    \end{aligned} $
    \;\vskip -12pt
    $ \pi_{i,i} \coloneqq \beta_i, \ \pi_{i,i-1} \coloneqq 1 $
    \;
    \KwSty{for} $k=i+1\ \KwSty{to}\ n$
    \KwSty{do} 
        $ \displaystyle \pi_{i,k} \coloneqq \pi_{i,k-1} \times \beta_k $
    \;
}
\For{j=1 \KwSty{to} n}{
    \KwSty{for} $k = j+2$ \KwSty{to} $n$ \KwSty{do} $ \displaystyle q_{j,j} \coloneqq q_{j,j} + \frac{ \pi_{j+2,k-1}^2}{\dot v_{k} \dot v_{k-1}}$
    \;
    \KwSty{if} $j < n - 1$ \KwSty{then} $ \displaystyle q_{j,j} \coloneqq q_{j,j} \times \left(\frac{\dot v_{j+1}-\beta_{j+1} \dot v_{j}}{b_{j+1}-b_j} - \beta_{j+1}\frac{\dot v_j-\beta_j \dot v_{j-1}}{b_j-b_{j-1}}\right)^2 $
    \;
    $ q_{j,j} \coloneqq q_{j,j} + \frac1{(b_j-b_{j-1})^2}\frac{\dot v_{j-1}}{\dot v_j} $
    \;
    \KwSty{if} $j < n$ \KwSty{then} $ \displaystyle q_{j,j} \coloneqq q_{j,j} + \left(
                \frac{\dot v_j - \beta_j \dot v_{j-1}}{b_{j}-b_{j-1}}
                +
                \frac{\dot v_{j}}{b_{j+1} - b_{j}}
            \right)^2
            \frac{1}{\dot v_{j+1}\dot v_j} $
    \;
    \For{i=1 \KwSty{to} j-1}{
        \KwSty{for} $k = j+2$ \KwSty{to} $n$ \KwSty{do} $ \displaystyle q_{i,j} \coloneqq q_{i,j} + \frac{ \pi_{i+2,k-1}\pi_{j+2,k-1}}{\dot v_{k} \dot v_{k-1}}$
        \;
        \KwSty{if} $j < n - 1$ \KwSty{then} $ \displaystyle q_{i,j} \coloneqq q_{i,j} \times \left(\frac{\dot v_{i+1}-\beta_{i+1} \dot v_{i}}{b_{i+1}-b_i} - \beta_{i+1}\frac{\dot v_i-\beta_i \dot v_{i-1}}{b_i-b_{i-1}}\right)\left(\frac{\dot v_{j+1}-\beta_{j+1} \dot v_{j}}{b_{j+1}-b_j} - \beta_{j+1}\frac{\dot v_j-\beta_j \dot v_{j-1}}{b_j-b_{j-1}}\right) $
        \;
        \If{$j < n$}{ 
        $
        \displaystyle
            q_{i,j} \coloneqq q_{i,j} - \left(\frac{\dot v_{i+1}-\beta_{i+1} \dot v_{i}}{b_{i+1}-b_i} - \beta_{i+1}\frac{\dot v_i-\beta_i \dot v_{i-1}}{b_i-b_{i-1}}\right)
            \left(
                \frac{\dot v_{j} - \beta_{j}\dot v_{j-1}}{b_{j}-b_{j-1}}
                + \frac{\dot v_{j}}{b_{j+1} - b_{j}}
            \right)
            \frac{\pi_{i+2,j}}{\dot v_{j+1}\dot v_j}
        $ 
        }
        \eIf{$i<j-1$}{
            $\displaystyle q_{i,j} \coloneqq q_{i,j} + \frac{1}{b_j-b_{j-1}}\left(\frac{\dot v_{i+1}-\beta_{i+1} \dot v_{i}}{b_{i+1}-b_i} - \beta_{i+1}\frac{\dot v_i-\beta_i \dot v_{i-1}}{b_i-b_{i-1}}\right) \frac{\pi_{i+2,j-1}}{\dot v_j}$
            \;
        }(\tcp*[h]{Case i = j-1}){
            $ \displaystyle q_{i,j} \coloneqq q_{i,j}   -\frac{1}{b_j-b_{j-1}}\left(
                \frac{\dot v_{j-1} - \beta_{j-1}\dot v_{j-2}}{b_{j-1}-b_{j-2}}
                + \frac{\dot v_{j-1}}{b_j - b_{j-1}}
            \right) \frac{1}{\dot v_j}
            $
        }
        %$ q_{j,i} \coloneqq q_{i,j} $ \;
    }
}
\caption{Peusdocode for the inverse of $\vec{A+C}$ using equations \eqref{eq:suminv-diagonal} and \eqref{eq:suminv-upper}.}
\end{algorithm}

\subsection{Numerical stability for a family of matrices with determinant $-\varepsilon/9$}

\begin{figure}[h]
    \centering
    \includegraphics{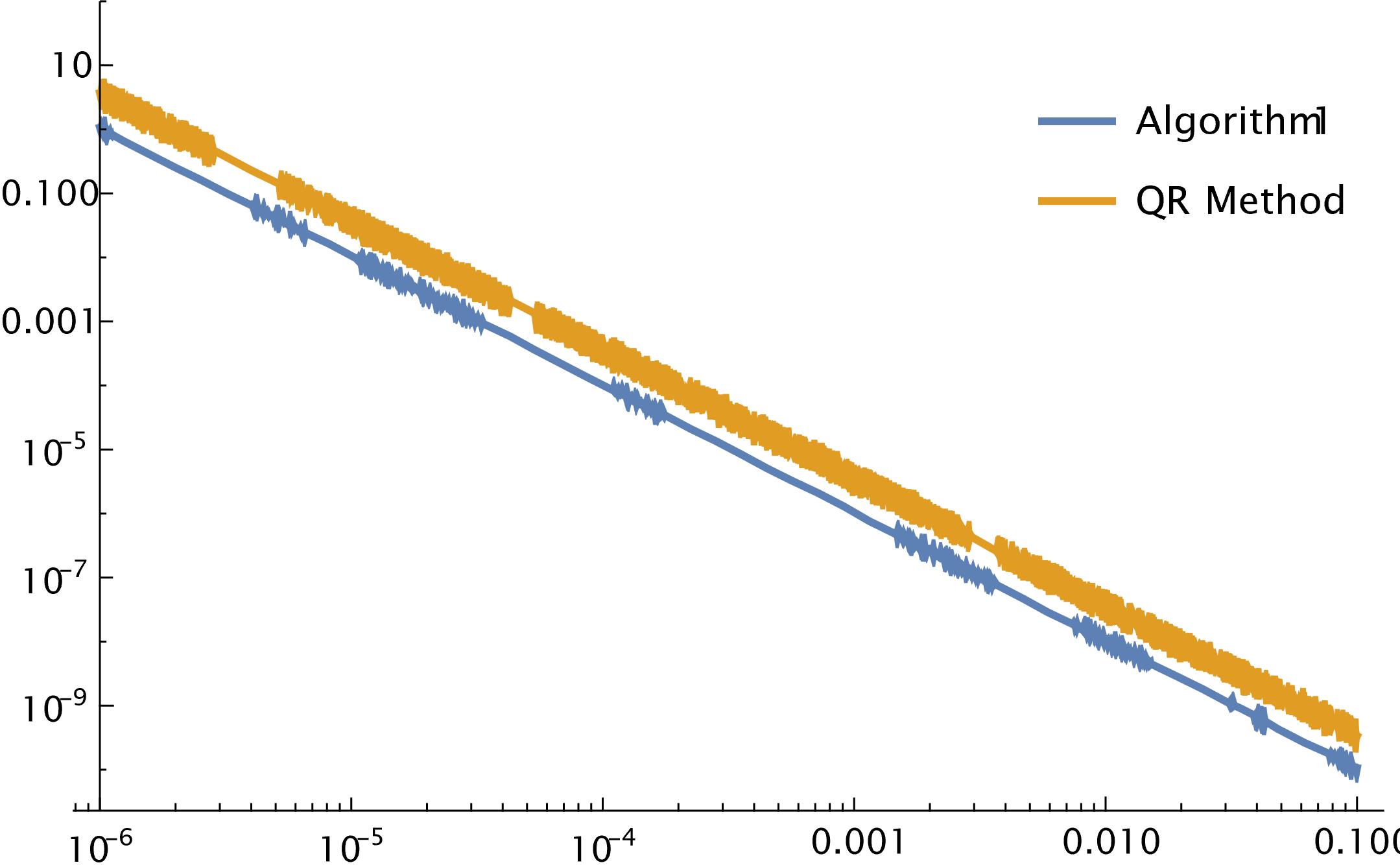}
    \caption{MAE of matrices with determinant $-\varepsilon/9$ as a function of $\varepsilon\in[10^{-6},0.1]$}
    \label{fig:stability-1}
\end{figure}

We use $\vec a = (1, 1, 1)^T, \vec b = (1,5/3,3)^T, \vec c(\varepsilon) = (0,1,\varepsilon-3)^T, x = 0, z = 1$, corresponding to
\[
    \vec{A+C(\varepsilon)} = \begin{pmatrix}
         1 & 5/3 & 3 \\
         5/3 & 8/3 & 4 \\
         3 & 4 & \varepsilon        
    \end{pmatrix},
    \qquad 
    (\vec{A+C(\varepsilon)})^{-1} = \begin{pmatrix}
        144/\varepsilon-24 & 15-108/\varepsilon & 12/\varepsilon \\
        15-108/\varepsilon & 81/\varepsilon-9 & -9/\varepsilon \\
        12/\varepsilon & -9/\varepsilon & 1/\varepsilon
    \end{pmatrix}.
\]
For each inverse matrix $\vec M(\varepsilon)$ computed by the algorithm, we calculate the mean absolute error (MAE) of coefficients $\vec M_{i,j}(\varepsilon)$ against the closed-form coefficients above, i.e.
\[
    \operatorname{MAE}(\varepsilon) = \frac19\sum_{i,j=1}^3 \norm*{\vec M_{i,j}(\varepsilon) - (\vec{A+C(\varepsilon)})^{-1}_{i,j}}.
\]
Figure \ref{fig:stability-1} reports the MAE as a function of $\varepsilon\in[10^{-6},0.1]$ on a log-log plot, together with the MAE of the QR method (i.e. the inverse $\vec M(\varepsilon)$ of $\vec{A+C}(\varepsilon)$ computed via QR decomposition).  We can see that for both inversion methods, the MAE exhibits a linear behavior with small and rapid oscillations that are typical of numerical rounding.  As expected, higher error is observed when $\varepsilon$ is lower, i.e. when the condition number is high.  In addition, our algorithm appears slightly more accurate for this particular family of matrices.

\subsection{Numerical stability for a class of matrices with spectrum $\{-1,\varepsilon, 1\}$}

In this approach we consider $\vec{A+C}$ with fixed generator $\vec b = (1,-1,1)^T$, and systematic generators $\vec a(\varepsilon), \vec c(\varepsilon)$ such that $\vec A(\varepsilon)+\vec C(\varepsilon)$ has eigenvalues $-1, \varepsilon, 1$.  Specifically, coefficients $a_1, a_3, c_2$ run from -1 to 1 with step 0.01 for a total of $201^3 \approx 8,120,000$ combinations, while $a_2, c_1, c_3$ are numerically solved for up to 6 real solutions so that $\operatorname{tr}(\vec{A+C}) = \varepsilon, \det(\vec{A+C})=-\varepsilon, \operatorname{tr}(\vec{A+C})^2 = 2 + \varepsilon^2 $.  Due to its computational cost, this method was only applied for three values of $\varepsilon = 10^{-1}, 10^{-4}, 10^{-7}$, resulting in approximately 27,000,000 real matrices.

\begin{figure}[p]
    \centering
    \includegraphics[width=0.65\textwidth]{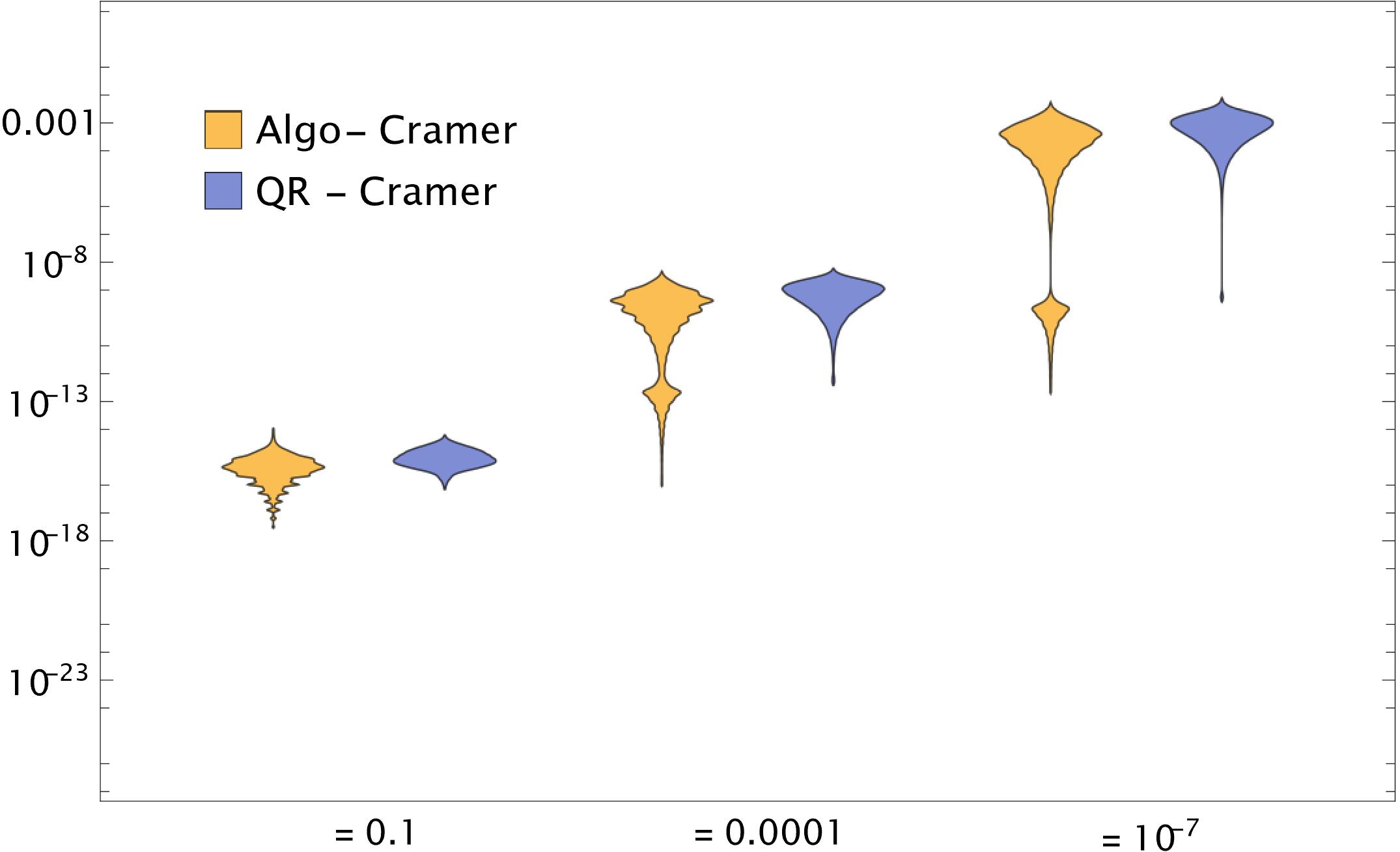}
    \caption{MAE distribution charts,  27mn. matrices with spectrum $\{-1, \varepsilon, 1\}$.}
    \label{fig:stability-2}
\end{figure}

\begin{figure}[p]
    \centering
    \includegraphics[width=0.65\textwidth]{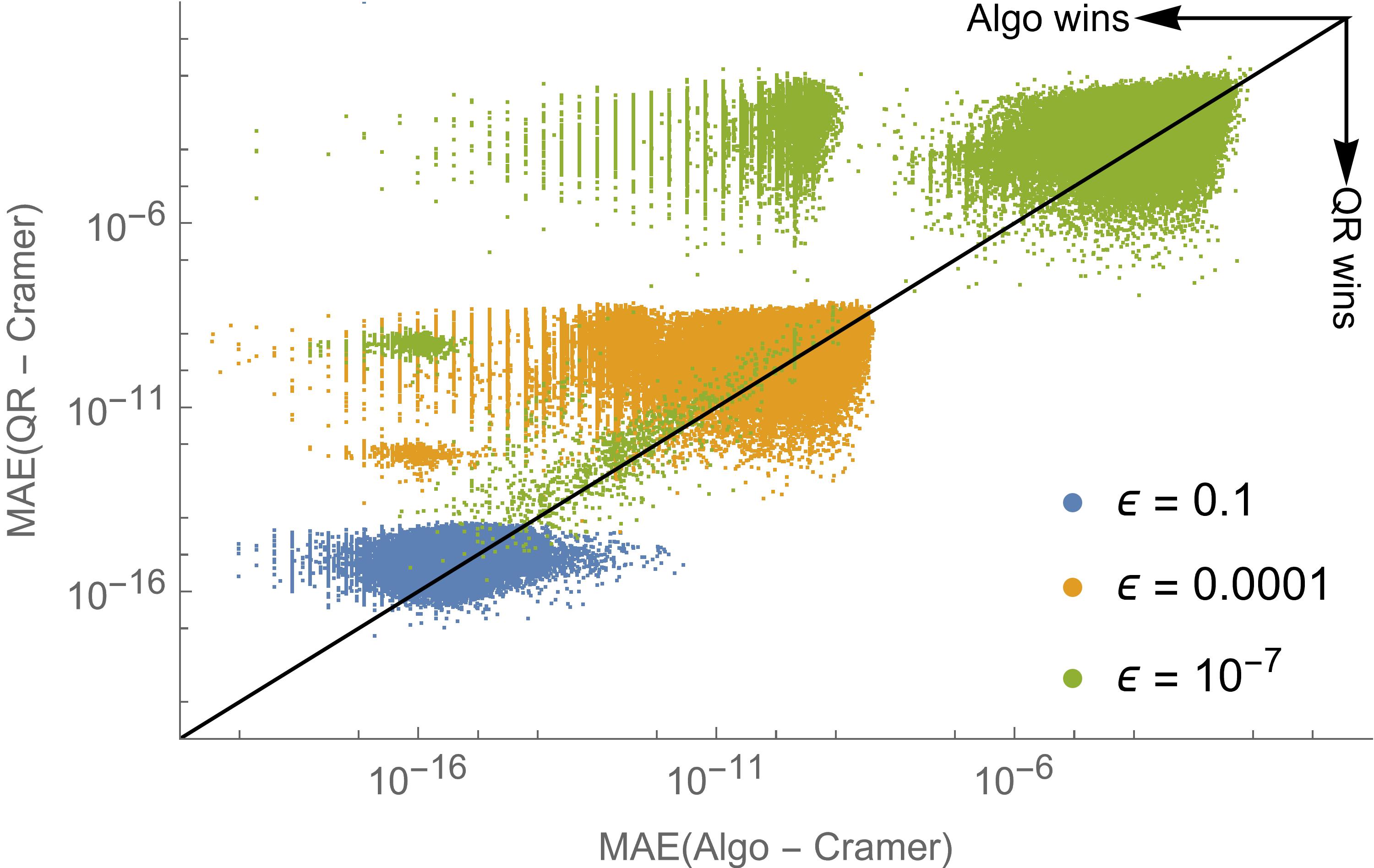}
    \caption{MAE scatter,  27mn. matrices with spectrum $\{-1, \varepsilon, 1\}$.}
    \label{fig:stability-3}
\end{figure}

For each matrix in our sample, we calculated the MAE of Algorithm 1 and the QR method versus inversion by Cramer's rule (i.e. the adjugate matrix divided by the determinant -$\varepsilon$). Table \ref{tab:stability-2} reports the corresponding statistics grouped by $\varepsilon = 10^{-1}, 10^{-4}, 10^{-7}$, while Figures \ref{fig:stability-2} and \ref{fig:stability-3} respectively show the distribution charts and scatter plots for each $\varepsilon$.  We can see that both methods perform similarly compared to Cramer's rule, with slightly better results for our algorithm.  However, while Algorithm 1 is often more accurate than QR, the wider distribution range may indicate less predictability in output, which could result in lower numerical stability for applications.

Overall, the above empirical studies suggest that Algorithm 1 is competitive versus the QR method.

\begin{table}[htb]
    \centering
    \begin{tabular}{|l|LL|LL|LL|}
        \hline
        & \multicolumn{2}{c|}{$\varepsilon =.1$} & \multicolumn{2}{c|}{$\varepsilon =.0001$} & \multicolumn{2}{c|}{$\varepsilon = 10^{-7}$}\Tstrut \\
        MAE & \text{Algo--Cramer} & \text{QR--Cramer} & \text{Algo--Cramer} & \text{QR--Cramer} & \text{Algo--Cramer} & \text{QR--Cramer}\Bstrut \\
        \hline
        Sample size& \multicolumn{2}{c|}{8,730,234} & \multicolumn{2}{c|}{9,143,088} & \multicolumn{2}{c|}{9,145,271}\Tstrut \\
        Average & 8.35\times10^{-16} & 1.08\times10^{-15} & 3.12\times10^{-10} & 8.62\times10^{-10} & 3.23\times10^{-4} & 8.90\times10^{-4}
        \\
        Std. dev. & 1.35\times10^{-14} & 8.26\times10^{-16} & 4.51\times10^{-10} & 8.16\times10^{-10} & 4.91\times10^{-4} & 9.30\times10^{-4} \\
        Minimum & 0. & 0. & 0. & 3.95\times10^{-15} & 0. & 2.04\times10^{-16} \\
        Median & 3.82\times10^{-16} & 8.26\times10^{-16} & 1.41\times10^{-10} & 6.26\times10^{-10} & 1.31\times10^{-4} & 6.11\times10^{-4}
        \\
        99\textsuperscript{th} pctile & 4.06\times10^{-15} & 3.94\times10^{-15} & 2.09\times10^{-9} & 3.57\times10^{-9} & 2.29\times10^{-3} & 4.13\times10^{-3} \\
        Maximum & 3.40\times10^{-12} & 8.68\times10^{-15} & 5.70\times10^{-9} & 8.56\times10^{-9} & 2.44\times10^{-2} & 5.70\times10^{-2}
        \Bstrut\\
        \hline
    \end{tabular}
    \caption{MAE sample statistics, 27mn. matrices with spectrum $\{-1, \varepsilon, 1\}$. }
    \label{tab:stability-2}
\end{table}


\clearpage
\printbibliography

\appendix
\section*{Appendix: implementation of Algorithm 1 in C code.}

\noindent\footnotesize{\tt{// License:\ CC BY-NC 4.0:\ Attribution-NonCommercial 4.0 Int'l \href{https://creativecommons.org/licenses/by-nc/4.0/}{creativecommons.org/licenses/by-nc/4.0}}}
\begin{lstlisting}[language=C,basicstyle=\ttfamily\footnotesize]
static double square(const double h) { return h * h; }
#define ABORT(exitcode, msg, ...){ fprintf(stderr, "Error %d: " #msg "\n", exitcode, __VA_ARGS__); \
                                    return exitcode; }
/// Calculates the inverse sum of two single-pair matrices A + C
/// n: Dimension
/// tol: Tolerance for zero tests
/// oa: Generator a[1]...a[n] of A (array of size n+1)
/// xb: Parameter b[0] = x and generator b[1]...b[n] of A (array of size n+1)
/// oc: Generator c[1]...c[n] of C (array of size n+1)
/// Q: Output upper-triangular inverse q(i,j) = Q[(2*n-i)*(i-1)/2+j-1] (array of size n(n+1)/2)
/// Pi: Output products of beta coefficients (array of size n(n+1)/2-1)
/// zv: Parameter v[0] = z, output continuants v[1]...v[n] (array of size n+1)
/// obeta: Output beta coefficients (array of size n)
/// verbose: Flag for warning output
/// Returns: 0 if successful, <0 if error, >0 if warning
#define q(i,j) Q[(2*n-(i))*(i-1)/2+j-1]
#define pi(i,j) Pi[(2*n+1-(i))*(i-2)/2+j-1]
#define a oa
#define b xb
#define c oc
#define v zv
#define beta obeta
int c_SPSumInverse(int n, double tol, double oa[], double xb[], double oc[],
                   double Q[], double Pi[], double zv[], double obeta[], _Bool verbose)
{
static int i, j ,k, warning = 0;
//Check arguments
if (!(a && b && c && Q && Pi && v && beta))
  ABORT(-1,"An array pointer is NULL");
if (n < 3 || !(isfinite(tol) && isfinite(b[0]) && isfinite(v[0])) || tol <= 0. || fabs(v[0]) < tol)
  ABORT(-2,"Invalid n = %d,tol = %e,x = b[0] = %e,or z = v[0] = %e",n,tol,b[0],v[0]);
if (fabs(b[1]-b[0]) < tol)
  ABORT(-3,"Invalid x = b[0] = %e too close to b[1] = %e",b[0],b[1]);
//Initialize
a[0] = c[0] = beta[0] = 0.;
v[1] = (v[0]*(a[1]*b[1]+c[1]))/square(b[1]-b[0]);
if (fabs(v[1]) < tol)
  ABORT(-4,"Invalid x = b[0],z = v[0],a[1],b[1],or c[1] (v[1] = %e ~ 0)",v[1]);
beta[1] = (a[1]*b[0]+c[1])/square(b[1]-b[0]);
//Compute continuants and betas
for (i = 2; i <= n; i++) {
  if (fabs(b[i]-b[i-1]) < tol || fabs(b[i-1]-b[i-2]) < tol)
    ABORT(-16,"Some consecutive b[] values are too close around index %d",i);
  v[i] = ((a[i]*b[i]-2*a[i-1]*b[i]+a[i-1]*b[i-1]+c[i]-c[i-1])/square(b[i]-b[i-1])
         - (a[i-1]*b[i-1]-2*a[i-1]*b[i-2]+a[i-2]*b[i-2]-(c[i-1]-c[i-2]))/square(b[i-1]-b[i-2]))*v[i-1]
         - square((a[i-1]*b[i-2]-a[i-2]*b[i-1]+c[i-1]-c[i-2]) / square(b[i-1]-b[i-2]))*v[i-2];
  if (fabs(v[i]) < tol)
    ABORT(-17,"Low continuant error (v[%d] = %e ~0)",i,v[i]);
  beta[i] = (a[i]*b[i-1]-a[i-1]*b[i]-c[i-1]+c[i])/square(b[i]-b[i-1]);
  if (fabs(v[i]-beta[i]*v[i-1]) < tol) {
    if (verbose) fprintf(stderr,"Warning: continuant value v[%d] ~ beta[%d]v[%d] might cause \
unreliable results.\n",i,i,i-1);
    warning = 18;
  }
  pi(i,i) = beta[i];
  pi(i,i-1) = 1.;
}
  //Compute products of betas
for (i = 2; i <= n; i++)
  for (k = i+1; k <= n; k++)
    pi(i,k) = pi(i,k-1)*beta[k];
  //Compute inverse
for (j = 1; j <= n; j++) {
  q(j,j) = 0.;
  if (j < n-1) {
    for (k = j+2; k <= n; k++)
      q(j,j) += square(pi(j+2,k-1))/(v[k]*v[k-1]);
    q(j,j) *= square((v[j+1]-beta[j+1]*v[j])/(b[j+1]-b[j])-(beta[j+1]*(v[j]-beta[j]*v[j-1]))
                      /(b[j]-b[j-1]));
  }
  q(j,j) += v[j-1]/(square(b[j]-b[j-1])*v[j]);
  if (j < n)
    q(j,j) += square((v[j]-beta[j]*v[j-1])/(b[j]-b[j-1])+v[j]/(b[j+1]-b[j]))/(v[j+1]*v[j]);
  for (i = 1; i <= j-1; i++) {
    q(i,j) = 0.;
    if (j < n-1) {
      for (k = j+2; k <= n; k++)
        q(i,j) += (pi(i+2,k-1)*pi(j+2,k-1))/(v[k]*v[k-1]);
      q(i,j) *= ((v[i+1]-beta[i+1]*v[i])/(b[i+1]-b[i]) - 
                    (beta[i+1]*(v[i]-beta[i]*v[i-1]))/(b[i]-b[i-1])) *
                    ((v[j+1]-beta[j+1]*v[j])/(b[j+1]-b[j])-(beta[j+1]*(v[j]-beta[j]*v[j-1]))
                    /(b[j]-b[j-1]));
    }
    if (j < n)
      q(i,j) -= (((v[i+1]-beta[i+1]*v[i])/(b[i+1]-b[i])-(beta[i+1]*(v[i]-beta[i]*v[i-1]))
                  /(b[i]-b[i-1]))*((v[j]-beta[j]*v[j-1])/(b[j]-b[j-1])+v[j]/(b[j+1]-b[j]))
                  *pi(i+2,j))/(v[j+1]*v[j]);
    if (i < j-1)
      q(i,j) += (((v[i+1]-beta[i+1]*v[i])/(b[i+1]-b[i])-(beta[i+1]*(v[i]-beta[i]*v[i-1]))
                  /(b[i]-b[i-1]))*pi(i+2,j-1))/((b[j]-b[j-1])*v[j]);
    else // i = j - 1
      q(i,j) -= ((v[j-1]-beta[j-1]*v[j-2])/(b[j-1]-b[j-2])+v[j-1]/(b[j]-b[j-1]))
                  /((b[j]-b[j-1])*v[j]);
  }
}
return warning;
}
\end{lstlisting}

\end{document}